\documentclass[a4paper,12pt]{article}
\usepackage{amssymb,amsmath,amscd,amsfonts,amsthm,mathtools}
\usepackage{latexsym,graphpap,xcolor}
\usepackage[utf8]{inputenc}
\usepackage[english]{babel}
\usepackage{enumerate}

\usepackage[right=3cm,left=3cm,top=3cm,bottom=3cm]{geometry}

\newtheorem{theorem}{Theorem}[section]
\newtheorem{prop}[theorem]{Proposition}
\newtheorem{lemma}[theorem]{Lemma}
\newtheorem{coro}[theorem]{Corollary}

\newtheorem{defi}[theorem]{Definition}

\newtheorem{result}[theorem]{Result}
\newtheorem*{remark}{Remark}
\newenvironment{demo}{ \noindent 
\emph{\textbf{Proof:}}}{\hfill$\square$\\}

\newcommand{\RR}{\mathbb{R}}
\newcommand{\NN}{\mathbb{N}}

\newcommand{\ZZ}{\mathbb{Z}}

\newcommand{\QQ}{\mathbb{Q}}

\newcommand{\Gg}{\mathfrak{G}}

\newcommand{\Ac}{\mathcal{A}}

\newcommand{\Cc}{\mathcal{C}}

\newcommand{\Fc}{\mathcal{F}}
\newcommand{\Gc}{\mathcal{G}}

\newcommand{\Lc}{\mathcal{L}}
\newcommand{\Mc}{\mathcal{M}}
\newcommand{\Nc}{\mathcal{N}}
\newcommand{\Oc}{\mathcal{O}}
\newcommand{\Pc}{\mathcal{P}}
\newcommand{\Qc}{\mathcal{Q}}

\newcommand{\Vc}{\mathcal{V}}

\renewcommand{\d}{\,{\rm d}}
\newcommand{\id}{\,{\rm id}}

\newcommand{\no}{n$^{\text{o}}$}

\newcommand{\Ker}{\operatorname{Ker}}
\newcommand{\R}{\operatorname{R}}

\newcommand{\Ind}{\operatorname{Ind}}
\newcommand{\pr}{\operatorname{p}}

\newcommand{\pc}{ \usefont{T1}{cmtl}{m}{n} \selectfont}

\newdimen\texpscorrection
\newdimen\figcenter
\def\figurewithtex #1 #2 #3 #4 #5\cr{\null
  {\goodbreak\figcenter=\hsize\relax
  \advance\figcenter by -#4truecm
  \divide\figcenter by 2
  \begin{figure}[hbt]
  \vskip #3truecm\noindent\hskip\figcenter
  \includegraphics{#1}{\hskip\texpscorrection\input #2 }
  \vskip 0.8truecm{\baselineskip=0.8\baselineskip
  \noindent \vbox{\noindent {\footnotesize #5}}\par}
  \end{figure}}}
\def\point#1 #2 #3 {\rlap{\kern #1 truecm
\raise #2 truecm \hbox{#3}}}

\numberwithin{equation}{section}

\begin{document}

\title{\bf Generic balanced synchrony patterns in network dynamics}

\author{Romain {\sc Joly} \& Maxime {\sc Percie du Sert}}

\maketitle

\bigskip

\begin{abstract}
{
A coupled cell network is a type of ordinary differential equation $\dot x(t)=f(x(t))$, with structural constraints on the vector field $f$, encoded in a directed graph, whose cells and arrows are labeled by type. The generated dynamics can model, for example, those of neural networks or ecological systems. These systems and the synchrony patterns observed in their solutions have been intensely studied, particularly by Golubitsky, Stewart, and their coauthors.

In the present article, we show that, for a generic vector field $f$, the synchrony patterns of the solutions of $\dot x(t)=f(x(t))$ are always balanced.  This roughly means that for almost all $f$, the observed synchrony patterns, such as synchronization in two different cells, are inherited from the structural symmetries imposed by the graph and the cell types. Any other synchronization, not directly imposed by the geometry of the graph and the cell types, cannot occur.

By doing so, we are completing the proof of several conjectures, including the rigid synchrony conjecture, the full oscillation conjecture and the observation of constant states.

This article is the published version of the results stated by the second author in his PhD thesis.
\\[5mm]}

\end{abstract}

\pagebreak

\section{Introduction}

\begin{figure}[ht]
\begin{center}
\resizebox{5cm}{!}{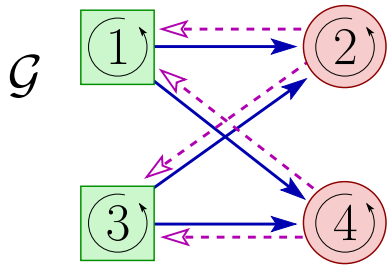} 
\end{center}
\caption{\it The graph $\Gc$ above has $4$ cells linked with $8$ arrows. There are two types of cells: the cells $1$ and $3$ (the left/green/squared ones) and the cells $2$ and $4$ (the right/red/round ones). There are two types of arrows: the ones from left to right (the blue/solid ones) and the ones from right to left (the magenta/dashed ones). Note that we include circling arrows inside the cells to remember that the evolution of a state also depends on itself.} \label{fig_exemple1}
\end{figure}

This introductory section is intended for readers who are not familiar with coupled cell networks and the associated dynamics. We will focus on a representative example provided by the graph of Figure \ref{fig_exemple1} and we will keep the discussion simple at the cost of being imprecise. 
In particular, our results will be stated in a fairly vague way and we will not reward them of the status of ``theorem''. The exact definitions used for our model and our main result are provided in Section \ref{section_notations} and we refer readers seeking greater rigor to this section. Other results will also be stated more precisely throughout the text. Readers familiar with network dynamics may jump directly to the statements of Theorem \ref{th_main}, Propositions \ref{prop_synchrony_eq} and  \ref{prop_rigid_constant} and Corollaries \ref{coro_rigid}, \ref{coro_phase_shift_1} and \ref{coro_phase_shift_2}. 

\medskip 

\noindent {\bf $\rhd$ Network dynamics}\\
Let $X=\RR^d$ with $d\geq 1$. Coupled cells networks are particular models of ODEs 
\begin{equation}\label{ode_intro}
\dot x(t)=f(x(t))
\end{equation}
where the structure of the vector field $f\in\Cc^k(X,X)$, $k\geq 1$, is constrained by the structure of a given directed graph $\Gc$ and the labeling of its cells and arrows by types:
\begin{enumerate}[(i)]
\item The state space $X$ is split between the cells of $\Gc$, that is that $X=\prod_{c} X_c$ where $X_c$ is the state space in the cell $c$. Furthermore, the directed arrows of $\Gc$ describe the possible inputs of a cell: if there is no arrow from a cell $c'$ pointing to a cell $c$ then the component $f_c$ of the vector field that governs the evolution of the state in $c$ must be independent on the state in $c'$.
\item If two cells $c$ and $c'$ have the same type and if their inputs are given by arrows and cells of the same type, then the components $f_c$ and $f_{c'}$ of the vector field must be the same. Note that this constraint also applies when $c=c'$: if, for example, a cell has two inputs of the same type, then the vector field must be invariant under the permutation of these two inputs. 
\end{enumerate}
To illustrate the above conditions, consider the example of the graph of Figure \ref{fig_exemple1}. For the sake of simplicity, assume that the state space of each cell is $\RR$. A state is of the form $x=(x_1,x_2,x_3,x_4)\in\RR^4$ where the component $x_i$ lives in the cell $i$. The graph $\Gc$ of Figure \ref{fig_exemple1} codes for the ODE \eqref{ode_intro} with the following constraints:
\begin{enumerate}[(a)]
\item there exist two functions $g,h\in\Cc^k(\RR^3,\RR)$ such that 
\begin{equation}\label{exemple_intro_1}
\left\{\begin{array}{l}
\dot x_1(t)=f_1(x(t))=g(x_1,x_2,x_4)\\
\dot x_2(t)=f_2(x(t))=h(x_2,x_1,x_3)\\
\dot x_3(t)=f_3(x(t))=g(x_3,x_2,x_4)\\
\dot x_4(t)=f_4(x(t))=h(x_4,x_1,x_3)
\end{array}\right.
\end{equation}
To illustrate the link with $\Gc$, note for example that the first component $f_1(x)$ of the vector field is of the form $g(x_1,x_2,x_4)$ and does not depend on $x_3$. This is due to the absence of the arrow $3\rightarrow 1$ in $\Gc$. Also note that the cells $1$ and $3$ are of the same type with the same types of input arrows. This explains why the third component $f_3$ is of the form $g(x_3,x_2,x_4)$ with the same function $g$ as the one defining $f_1$.
\item the types of the inputs of the cell $1$ are invariant under the exchange of the arrows and the same symmetry holds for the other cells. This symmetry translates into constraints on the vector field which are
\begin{equation}\label{exemple_intro_2}
\forall (\alpha,\beta,\gamma)\in\RR^3~,~~~g(\alpha,\beta,\gamma)=g(\alpha,\gamma,\beta)~~~\text{ and }~~~h(\alpha,\beta,\gamma)=h(\alpha,\gamma,\beta).
\end{equation}
\end{enumerate}

\medskip 

\begin{figure}[p!]
\begin{center}
\resizebox{0.99\textwidth}{!}{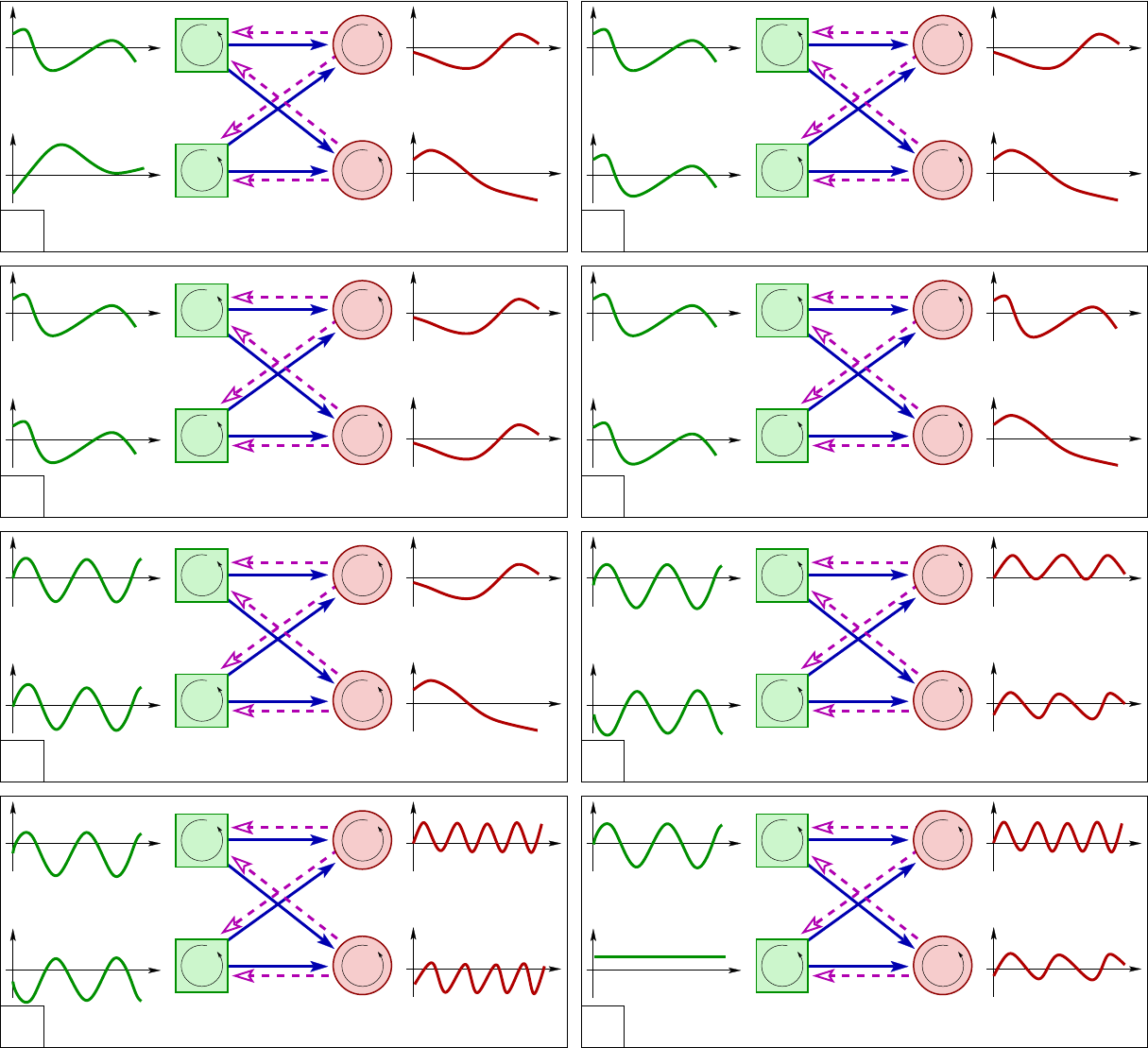} 
\end{center}
\caption{\it Several examples of the possible qualitative dynamics of the solutions of the network of Figure \ref{fig_exemple1}. The symmetries of the network—that is, of the graph $\Gc$ and its associated types—allow for solutions that present certain patterns. The purpose of the present article is to show that, generically with respect to the vector field, the non-expected patterns cannot appear. In the above figure, the cases A, B, C and G are compatible with the symmetries of the network and may be observed. Note that this list is incomplete: there are other balanced synchrony patterns that exist. On the opposite, the cases D, E, F and H are generically not possible because these synchrony patterns are not balanced, i.e. not compatible with the structural constraints and symmetries provided by the networks.}\label{fig_exemple2}
\end{figure}

\noindent {\bf $\rhd$ Symmetries and synchrony}\\
Consider a solution $t\mapsto x(t)$ of the ODE \eqref{exemple_intro_1}. If, for two cells $i$ and $j$, we have $x_i(t)=x_j(t)$ for all $t$, then we say that the cells $i$ and $j$ are synchronous and we write $i\bowtie j$. The relation $\bowtie$ is what we call a \emph{synchrony pattern}. 
Due to the constrained structure of the vector field, some of these synchrony patterns are expected to appear. For example, if we consider a solution $x(t)$ of the ODE \eqref{exemple_intro_1} with $x_1(0)=x_3(0)$, then $x_1(t)=x_3(t)$ for all times. This synchrony comes from the fact that the cells $1$ and $3$ are of the same type and have the same inputs of the same types. So the pattern $1\bowtie 3$ (case B of Figure \ref{fig_exemple2}) is a natural synchrony dictated by the structural constraints provided by $\Gc$ and we say that it is a \emph{balanced} synchrony pattern (see Section \ref{section_notations} for a more precise definition). This is also the case of the trivial synchrony where $i\bowtie j$ only if $i=j$ (case A of Figure \ref{fig_exemple2}) or the synchrony coupling both the left cells with $1\bowtie 3$ and the right cells with $2\bowtie 4$ (case C of Figure \ref{fig_exemple2}). On the contrary, a solution satisfying $x_1(t)=x_2(t)$ for all times has a priori no reason to exist because the vector field in the cell $1$ is not related to the one in the cell $2$. So, while the existence of solutions exhibiting the synchrony pattern $1\bowtie 2\bowtie 3$ (case D of Figure \ref{fig_exemple2}) is possible for a specific $f$, it is not expected for a general vector field.

This intuition is precisely the subject of our main result. 
\begin{result}{\rm The synchrony patterns are balanced}

Fix a graph $\Gc$ and associated types. For a generic vector field $f$ satisfying the constraints provided by the graph and the types, the only possible synchrony patterns of a solution of $\dot x(t)=f(x(t))$ are balanced. 

In the case of the graph of Figure \ref{fig_exemple1}, for a generic set of functions $g$ and $h$, the only possible synchrony patterns of a solution $x(t)$ of \eqref{exemple_intro_1} are: no synchrony, the synchrony $x_1(t)=x_3(t)$, the synchrony $x_2(t)=x_4(t)$ and the synchrony $(x_1(t)=x_3(t)) \wedge (x_2(t)=x_4(t))$.
\end{result}

This result is stated more precisely and rigorously in Theorem \ref{th_main} and is proven in Section \ref{sect_traj}. In short, this main result says that observing the synchrony of a trajectory provides information about the network's structure. 
The literature provides many examples of interesting synchrony, see \cite{GNS,GS_vs,GS_book,Wang} for examples. To illustrate the practical interest of the above result, consider a coupled cell network that models a neural network, the graph $\Gc$ describing the neurons connections. Since the interaction between neurons is not expected to be governed by laws as precise as those of physics, it is reasonable to expect that the associated vector field $f$ is generic among all the vector fields that satisfy the constraints required by $\Gc$. Suppose that we are observing a group of neurons that have synchronous behaviour, for example that fire simultaneously. The above result provides a strong hint that we should consider that these neurons are similar (because modeled by cells of the same type) and have similar connections (because modeled by cells with the same arrows in $\Gc$).

\medskip 

\noindent {\bf $\rhd$ Oscillations and phase-shift}\\
Several observation properties can be deduced for the above result. Consider for example the following question: if the state $x_1(\cdot)$ in cell $1$ is constant, is the whole state $x(\cdot)$ constant in all cells? In \cite{RJ}, the first author provides a positive answer for generic networks without type constraints (in the framework of the present paper, this means assuming that all cells and arrows are of different types). 
More precisely, the property of being constant is expected to propagate upstream: if a solution $x(\cdot)$ is constant in a cell $i$, it should also be constant in any cell $j$ being one of its input, that is a cell for which there is an arrow $j\rightarrow i$. Consequently, it propagates to all indirect inputs and, if the graph is transitive, to the entire graph. 

In Proposition \ref{prop_rigid_constant} below, we extend this fact to the case of networks with symmetry. We can call this property an "observation" of stationary states, in analogy with the observation properties in control theory. Indeed, from partial information—the observation of the solution in a single cell—we recover global information: the qualitative behavior of the solution in all cells.   
\begin{result}{\rm Observation of stationary states}

Fix a graph $\Gc$ and associated types. For a generic vector field $f$ satisfying the constraints provided by the graph and the types, if a solution of $\dot x(t)=f(x(t))$ is constant in a cell $i$, then $x_j(\cdot)$ is also constant for any cell $j$ which is an indirect input of $i$.

In the case of the graph of Figure \ref{fig_exemple1}, for a generic set of functions $g$ and $h$, 
and for any solution $x(t)$ of \eqref{exemple_intro_1}, if there is an interval $J$ and a cell $i$ such that $t\in J\mapsto x_i(t)$ is constant, then the whole solution $t\in\RR\mapsto x(t)$ is constant.
\end{result}
In particular, it is generically impossible to have a non-constant periodic solution $t\mapsto x(t)$ which is constant in a cell $i$: the case H of Figure \ref{fig_exemple2} is generically impossible. The previous property is called the \emph{full oscillation property}, but we can see that it is actually a consequence of a more general property where the periodicity is not required. Also note that we cannot propagate the property of being constant downstream: for a two-cells graph as $1\rightarrow 2$, the state $x_1(\cdot)$ can be constant without precluding the motion of the state $x_2(\cdot)$. 

As for stationary states, the periodicity of the state in a cell propagates to the indirect input cells, as shown in Corollary \ref{coro_phase_shift_2} below.
\begin{result}{\rm Observation of periodic states}

Fix a graph $\Gc$ and associated types. For a generic vector field $f$ satisfying the constraints provided by the graph and the types, if a solution of $\dot x(t)=f(x(t))$ is periodic in a cell $i$, then $x_j(\cdot)$ is also periodic for any cell $j$ which is an indirect input of $i$.

In the case of the graph of Figure \ref{fig_exemple1}, for a generic set of functions $g$ and $h$, if a solution $x(t)$ of \eqref{exemple_intro_1} is such that $t\mapsto x_i(t)$ is periodic in a single cell $i$, then the whole solution $t \mapsto x(t)$ is periodic.
\end{result}
This means that a situation as case E of Figure \ref{fig_exemple2} is impossible for a generic vector field. 

Finally, we discuss another type of symmetry patterns that may occur in coupled cell networks, as shown in \cite{GS2} for example. 
Consider a periodic solution $x(\cdot)$ of minimal period $T$. Since the ODE $\dot x(t)=f(x(t))$ is autonomous, $x(\cdot+\theta)$ is also a periodic solution for any $\theta\in(0,T)$. A pattern of phase-shift synchrony appears when $x(\cdot+\theta)=\pi x(\cdot+\theta)$, where $\pi$ is a permutation of the cells. In the example of Figure \ref{fig_exemple1}, we can for example assume that $x_3(t)=x_1(t+T/2)$, as in the cases F and G of Figure \ref{fig_exemple2}. The \emph{balanced shift-phase property}  states that this phase shift is only possible if the permutation $\pi$ is compatible with the graph $\Gc$ and the types of cells and arrows: $x_{c'}(\cdot+\theta)=x_{c}(\cdot)$ only if both cells have the same type and the same type of inputs (this is the case for cells $1$ and $3$ of Figure \ref{fig_exemple1}) and only if we can find the same phase-shift in the inputs. This last property does not hold in case F of Figure \ref{fig_exemple2} since $x_2(\cdot)$ and $x_4(\cdot)$ are not related by a $T/2-$shift. But the case G of Figure \ref{fig_exemple2} is balanced: the inputs $x_2(\cdot)$ and $x_4(\cdot)$ of the cell $3$ are $T/2-$shifted from the inputs of the cell $1$ since both inputs are $T/2-$periodic. By the way, we notice that this situation is not exceptional: due to the $(\ZZ/2\ZZ)\times(\ZZ/2\ZZ)$ symmetry of \eqref{exemple_intro_2}, the inputs of the right cells are invariant by exchanging the left cells and so these left cells yield a $T/2-$periodic forcing for the right cells, which is compatible with their period. 

\begin{result}{\rm The phase-shift patterns are balanced}

Fix a graph $\Gc$ and associated types. For a generic vector field $f$ that satisfies the constraints provided by the graph and the types, the only possible phase-shift patterns of a periodic solution of $\dot x(t)=f(x(t))$ are balanced.

In the case of the graph of Figure \ref{fig_exemple1}, for a generic set of functions $g$ and $h$, if for example $x_3(t)=x_1(t+\theta)$, then $x_1(t)=x_3(t+\theta)$ and either $x_4(t)=x_2(t+\theta)=x_4(t+2\theta)$ or $x_2(t)=x_2(t+\theta)$ and $x_4(t)=x_4(t+\theta)$. In both cases, the whole trajectory is $2\theta-$periodic. 
\end{result}

\pagebreak

\medskip 

\noindent {\bf $\rhd$ Previous works}\\
The coupled cell networks are relevant models for many concrete dynamical systems, from chemical reactions \cite{Feinberg} or coupled oscillations \cite{BBH} to neural network \cite{ACN} and animal locomotion \cite{SW}. The study of the generic synchrony patterns is particularly motivated by the study of neural networks and coupled cell networks are very relevant toy models in this field. There are many publications on this subject.  For nice reviews of the literature, we refer to \cite{Stewart_nature} and the introductions of \cite{GS2} or \cite{Stewart}.

Before we discuss the literature, we must highlight a specific aspect of the timeline. The main results of the present article were stated and proven in the second author's PhD manuscript more than ten years ago, see \cite{Maxime}. However, the manuscript is written in French and has not been published in a research journal. Consequently, these results have remained ignored by the community and this explains why authors of recent articles still consider what was proven in \cite{Maxime} to be conjectures. This article aims to address this issue, and can be considered a revised version of some of the results of \cite{Maxime}, intended for publication in English. 
In the following discussions, we mostly adopt the perspective of the community and discuss the literature as if the results of the present article were new. 

The present article studies the synchrony patterns that the solutions of the differential system may satisfy, see \cite{Wang}, \cite{Stewart_nature} or \cite{GS_book} for several examples of applications.
In previous studies, this question has often been limited to an analysis of rigid patterns, which are synchronous patterns observed in a periodic orbit that are robust with respect to perturbations of the vector field, see Section \ref{sect_discuss_rigid}. The fact that the rigid synchrony patterns must be balanced has been conjectured in \cite{GS2,JT,SP}. The first proofs have been given in \cite{SP,GRW,GRW2} but they are only partial as underlined in the appendix of \cite{Stewart_eq}, see also \cite{Stewart}. Corrected proofs of the rigid synchrony conjecture are given in \cite{Stewart_eq} for equilibria and in \cite{Stewart} for periodic orbits (see also \cite{GS3}), but with the assumptions that these orbits are hyperbolic. Hyperbolicity is known to be generic within the entire class of all ODEs: this is part of the famous Kupka-Smale property, see for instance \cite{PM} or any handbook of dynamical systems. But the problem is that it is not known to be generic within the specific class of coupled cell networks. Thus, it is unclear whether the results of \cite{Stewart_eq,Stewart} are general or rely on too strong assumptions. However, we acknowledge that these previous works already contain almost all the ideas that we will use to tackle the difficulties coming from the geometry of coupled cell networks. As explained above, these ideas have been developed independently of \cite{Maxime}.

To be able to skip the assumption of hyperbolicity, we have to use the transversality theorems, also known as Sard-Smale theorems or Thom's theorems. These tools have been developed to prove generic results in geometry and in dynamical systems, see \cite{Abraham-Robbin,PM,Sard,Thom1,Thom2}. We will use here the version introduced by Henry in \cite{Henry}, see Section \ref{sect-trans}. In the previous work \cite{RJ}, the first author already used these techniques to prove the generic absence of synchrony in the case of fully inhomogeneous networks, that are networks with no imposed symmetries. The PhD thesis \cite{Maxime} and the present article extend the arguments to the case of networks with symmetries.  

We finish by enhancing that the genericity of the Kupka-Smale property has been proved in the case of fully inhomogeneous networks in \cite{Maxime}. It is noteworthy that \cite{Maxime} uses the results of \cite{RJ}. Thus, we could say that the generic hyperbolicity of periodic orbits should be obtained as a consequence of the generic balanced synchrony rather than being an assumption for proving it. Following \cite{Maxime} and the present article, we could hope to prove the genericity of Kupka-Smale property in the class of coupled cell networks with symmetries.

\vspace{5mm}

\noindent {\bf Acknowledgement:} We thank the anonymous referee for his careful reading and his numerous comments, which have helped us to improve our manuscript.


\section{The network dynamics: definitions, notations and main result}\label{section_notations}

{
This section introduces the main concepts and notations used in the present article. Most of these concepts and notations are classical and come from the series of works by Golubitsky, Stewart, and their co-authors, see \cite{GST,SGP} or the book \cite{GS_book}. Note that other formalisms exist, as the one of \cite{Field}, or the slightly more general one of \cite[Section 4]{Stewart} (see also \cite{GS_book}).}

\subsection{The network with types}\label{section_def_graph}

\begin{figure}[ht]
\begin{center}
\resizebox{0.9\textwidth}{!}{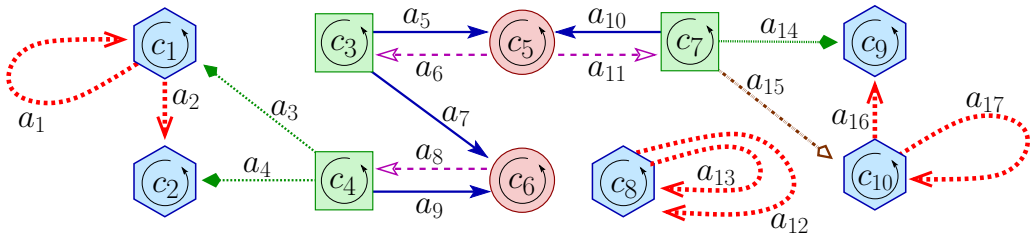} 
\end{center}
\caption{\it An example of network with types. The graph $\Gc$ has $10$ cells, divided in $3$ types, and $17$ arrows divided in $5$ types. The types are coded by shapes and colors. The small circling arrows inside each cell is a reminder that an internal arrow $c_i\rightarrow c_i$ is implicitly present, having its own type, see Definition \ref{defi_vectfield} below.}\label{fig_exemple3}
\end{figure}

The basic geometry of our dynamical systems is defined by a directed graph. Since the intended applications involve neural networks and economic networks, among others, we prefer to call the graph a ``network''. 
\begin{defi}
A {\bf network} $\Gc$ is a directed graph. It consists in:
\begin{enumerate}[(i)]
\item a set $\Cc$ of {\bf cells},
\item a set $\Ac$ of {\bf arrows} and functions $H:\Ac\rightarrow \Cc$ and $T:\Ac\rightarrow \Cc$ providing the {\bf head} $H(a)$ and the {\bf tail} $T(a)$ of each arrow. 
\end{enumerate}
\end{defi}
Consider the network of Figure \ref{fig_exemple3}. The arrow $a_3$ connects the cells $c_4$ to the cell $c_1$. In our notations, this exactly means $H(a_3)=c_1$ and $T(a_3)=c_4$. Notice that the arrow $a_1$ is such that $T(a_1)=H(a_1)=c_1$, which is perfectly fine. It is also possible to have multiple arrows as $a_{12}$ and $a_{13}$, i.e. several arrows having all the same head and the same tail. Also note that the network of Figure \ref{fig_exemple3} is not connected, which is not a problem for our study.

The dynamics inside a cell $c$ will be determined by the state inside the input cells of $c$.   
\begin{defi}\label{defi_input}
If $c$ is a cell of a network $\Gc$, the {\bf inputs} of $c$ are the elements of the set $I(c)$ of all the arrows $a$ pointing at $c$, i.e. $I(c)=\{a\in\Ac, H(a)=c\}$. The {\bf input cells} of $c$ are the cells directly connected to $c$, i.e. the elements of $T(I(c))$. 
\end{defi}
For example, in the network of Figure \ref{fig_exemple3}, the inputs of $c_1$ are $a_1$ and $a_3$ and the input cells of $c_1$ are $c_1$ and $c_4$. 

\medskip 

To the structure of directed graph, we add types on the cells and arrows.
\begin{defi}\label{defi_types}
A {\bf network with types} is a network $\Gc$ whose cells and arrows are classified into types. More precisely, there exists two relations of equivalence:
\begin{enumerate}[(i)]
\item The set of cells $\Cc$ is endowed with a relation of equivalence $\sim_\Cc$. We say that two cells $c$ and $c'$ are of the same type if $c\sim_\Cc c'$. 
\item The set of arrows $\Ac$ is endowed with a relation of equivalence $\sim_\Ac$. We say that two arrows $a$ and $a'$ are of the same type if $a\sim_\Ac a'$. 
\item The following compatibility condition is assumed to hold. If $a\sim_\Ac a'$ then $H(a)\sim_\Cc H(a')$ and $T(a)\sim_\Cc T(a')$, that is that two arrows of the same type connect the same type of cells. 
\end{enumerate}
\end{defi}
In the network of Figure \ref{fig_exemple3}, the arrows $a_5$ and $a_{10}$ are of the same type and the condition (iii) of Definition \ref{defi_types} holds since their heads and tails are of the same type. The arrows $a_3$ and $a_{15}$ connect cells of the same type but the two arrows have different types and this is fine because the reciprocal condition of (iii) is not mandatory. Note that it is possible to consider a slightly relaxed version of (iii), see \cite[Remark 2.2]{Stewart_eq} and \cite[Section 4]{Stewart}. It would allow us to include additional networks in our framework. However, for the sake of simplicity, we choose to keep the more traditional condition (iii).

\begin{defi}\label{defi_input_iso}
Let $\Gc$ be a network with types and $c$ and $c'$ two (possibly equal) cells. An {\bf input isomorphism} from $c$ to $c'$ is a bijective function from $I(c)$ into $I(c')$ such that $a\sim_\Ac \beta(a)$ for all $a\in I(c)$, i.e. $\beta$ is preserving the type of the input arrows. The set of all input isomorphisms from $c$ to $c'$ is denoted $B(c,c')$. If $B(c,c')\neq \emptyset$, we say that $c$ and $c'$ are {\bf input isomorphic}, which implies in particular that $c$ and $c'$ are cells of the same type.
\end{defi}

Consider again the network of Figure \ref{fig_exemple3} as an example. The three cells $c_3$, $c_4$ and $c_7$ are of the same type, each have a unique input arrow and this arrow has the same type for each cell. Thus $c_3$, $c_4$ and $c_7$ are input isomorphic cells with the obvious isomorphisms associating their respective unique input. As all cells, the cell $c_5$ is input isomorphic to itself because of the identity isomorphism. More interestingly, there is also a non trivial input isomorphism: the permutation of the arrows $a_5$ and $a_{10}$. The cells $c_5$ and $c_6$ are also input isomorphic and $B(c_5,c_6)$ contains two isomorphisms: either $a_5\mapsto a_7$ and $a_{10}\mapsto a_9$ or $a_5\mapsto a_9$ and $a_{10}\mapsto a_7$. On the contrary, $c_9$ and $c_{10}$ are not input isomorphic: even if the cells have the same type and admit cells of the same type as input cells, the input arrows $a_{14}$ and $a_{15}$ have different types.

\subsection{The associated ODE with symmetries}\label{section_def_ODE}
Consider a network with type $\Gc$. We associate a state space to this network as follows. To each cell $c$ is associated a state space $X_c$ and, in this article, we assume that $X_c=\RR^{d_c}$ for some dimension $d_c\geq 1$ (see Section \ref{sect_discuss_manifolds} for a discussion about other possible state spaces). We set $X=\prod_{c\in\Cc} X_c$, that is $X=\RR^d$ with $d=\sum_{c\in\Cc} d_c$. If $C=\{c_1,c_2,\ldots\}$ is a set of cells, we write $X_C$ for the subspace $X=\prod_{c\in C} X_c$ and, for $x\in X$, $x_{C}$ denotes the projection of $x$ on $X_C$, written more shortly $x_c$ if $C=\{c\}$. If two cells $c$ and $c'$ are such that $d_c=d_{c'}$ is relevant to think $X_c$ and $X_{c'}$ as two different copies of $\RR^{d_c}$ in order to distinguish the states inside each cell. However, we will often abusively use expressions as $x_c=x_{c'}$, by implying its obvious meaning and omitting any canonical identification map. 
It is convenient to endow all the spaces $\RR^k$ with the supremum norm $\|x\|=\max_{i} |x_i|$. In particular, we can write $\|x\|=\max_{c\in\Cc} \|x_c\|$ for $x\in X$. 

We have just endowed the network of a state space $X$. If the network $\Gc$ is endowed with cell- and arrow-types, we translate these symmetries on $X$.
\begin{defi}
Let $\Gc$ be a network with types. We say that $X=\prod_{c\in\Cc} X_c$ is an {\bf admissible state space} if $X_c=X_{c'}$ for any equivalent cells $c\sim_\Cc c'$.
\end{defi}
In an admissible state space, we can compare two equivalent cells. If we want to compare the inputs of two cells, we need to assume that they are input isomorphic to ensure that the state spaces of the input cells are comparable. Due to the  preservation of types in Definition \ref{defi_input_iso} above, the following map is well defined.
\begin{defi}
Let $\Gc$ be a network with types and $X$ an admissible state space. Let $c$ and $c'$ be two input isomorphic cells and $\beta\in B(c,c')$ an input isomorphism. 
We define the {\bf pullback map} $\beta^*$ as follows: for any list $(a_1,\ldots,a_p)$ of input arrows of $c$, we set  
$$\beta^*(x_{T(a_1)},\ldots , x_{T(a_p)})=(x_{T(\beta(a_1))},\ldots , x_{T(\beta(a_p))}).$$
\end{defi}
Next, we consider vector fields defined on $X$. Since $X$ is of the form $\RR^d$, we identify $X$ and its tangent space and the vector fields $f$ are function defined from $X$ to $X$. For such a vector field, we use the notation  $f_c:X\rightarrow X_c$ for the component of the function in the cell $c$. 
\begin{defi}\label{defi_vectfield}
Let $\Gc$ be a network with types and $X=\prod_{c\in\Cc} X_c$ an admissible state space. Let $k\geq 1$, we denote by $\Cc^k_\Gc$ the set of {\bf admissible vector fields}, that are the functions of class $\Cc^k(\RR^d,\RR^d)$ such that:
\begin{enumerate}[(i)]
\item For every cell $c$, the component $f_c$ only depends on the values of $x$ in the cell $c$ and its input cells $T(I(c))$. In other words, we assume that there exists a function $\hat f_c:X_c\times X_{T(I(c))}\rightarrow X_c$ such that 
$$f_c(x)=\hat f_c(x_c,x_{T(I(c))}).$$
\item If $c$ and $c'$ are two input isomorphic cells, then for every input isomorphism $\beta\in B(c,c')$, we have 
\begin{equation}\label{hyp_f_symmetry}
\hat f_{c}(x_c,x_{T(I(c))}) = \hat f_{c'}(x_c,\beta^* x_{T(I(c))})
\end{equation}
where $\beta^*$ is the pullback map defined above.
\end{enumerate}
\end{defi}
For an admissible vector field $f$, we abusively write $f_c(x)=f_c(x_c,x_{T(I(c))})$, that is that we omit the hat in the notation. 
As we can see, we include the self-dependence of all the states: $f_c$ depends on $x_c$ and this dependence is free from any imposed symmetry. This choice is the classical one of the previous works as \cite{GRW, GST,  Stewart_eq,Stewart} and we choose to keep it. It is also helpful in some of our proofs, even if it may be not mandatory, see Section \ref{sect_discuss_self_depend} for a short discussion on this self-dependence. In the network of Figure \ref{fig_exemple3}, this self-dependence is shown by the small internal circling arrows. Let us again use this example to illustrate Definition \ref{defi_vectfield} and consider an admissible vector field $f$. Then, it must be of the form
$$f(x)=\left(\begin{array}{c} f_1(x_1,x_1,x_4)\\f_1(x_2,x_1,x_4)\\f_3(x_3,x_5)\\f_3(x_4,x_6)\\f_5(x_5,x_3,x_7)\\f_5(x_6,x_3,x_4)\\f_3(x_7,x_5)\\f_8(x_8)\\f_1(x_9,x_{10},x_7)\\f_{10}(x_{10},x_7)    \end{array}\right)~~~ \text{ with }f_5(\zeta,\xi,\xi')=f_5(\zeta,\xi',\xi).$$
Indeed, let $f_i$ be the component of $f$ in the cell $c_i$. We note for example that: 
\begin{itemize}
\item The only input cell of the cell $c_3$ is $c_5$, so $f_3$ is defined from $X_3\times X_5$ into $X_3$ (the self-dependence always coming first). 
\item The cells $c_3$, $c_4$ and $c_7$ are input isomorphic cells and thus $f$ must satisfy $f_3(\zeta,\xi)=f_4(\zeta,\xi)=f_7(\zeta,\xi)$. 
\item The cell $c_5$ admits a non trivial internal input isomorphism consisting in exchanging the arrows $a_5$ and $a_{10}$. So we must have the symmetry $f_5(\zeta,\xi,\xi')=f_5(\zeta,\xi',\xi)$.
\item  Formally, the dynamics of the cell $c_8$ must follow a vector field of the type $f_8(x)=\hat f_8(x_8,x_8,x_8)$ due to both arrows $a_{12}$ and $a_{13}$. As we can see, this kind of doubled arrow is simply formal in the present paper, but it is interesting in other studies, see \cite{GST} for example. 
\item On the contrary, note that the arrow $a_{1}$, which connects $c_1$ to itself is not redundant with the fact that $f_1$ already depend on $x_1$ by Assumption (i) of Definition \ref{defi_vectfield}. Indeed, the cells $c_1$ and $c_2$ are input equivalent by the association $a_1\mapsto a_{2}$ and $a_3\mapsto a_{4}$ and so we must have the symmetry $f_2(\zeta,\xi,\xi')=f_1(\zeta,\zeta,\xi')$.
\item Finally, note that $c_{10}$ has the same types of input cells as the cell $c_1$ but the arrows $a_3$ and $a_{15}$ are of different types. So the cells are not input isomorphic and $f_{10}$ can be chosen independently of $f_1$. 
\end{itemize}

\medskip 

We consider in the present paper the ODE of the type
\begin{equation}\label{ode}
\dot x(t)=f(x(t)) 
\end{equation}
where $f\in\Cc^k_\Gc$ is an admissible $\Cc^k-$map, with $k\geq 1$ of a given network with types $\Gc$. { Note that the above ODE is locally well-posed by Cauchy-Lipschitz theory: $f$ is of class at least $\Cc^1$ and is thus a locally Lipschitzian function}. Our purpose is to study the dynamics of the flow generated by \eqref{ode} and in particular to understand the relation between the symmetries required by the network with types and the symmetries of the solutions of \eqref{ode}.

\medskip 

Since we aim at proving results that are generic with respect to $f$, we need to endow $\Cc^k_\Gc$ with a topology. To this end, we first consider the space 
$\Cc^k_{b,\Gc}$ of admissible functions $f$ such that every derivative $\partial^\alpha f$ with $|\alpha|\leq k$ is bounded on $X$. { Here and in what follows, we use the multi-index notation: $\alpha$ belongs to $\NN^d$, $\partial^\alpha$ stands for $\partial_{x_1}^{\alpha_1}\partial_{x_2}^{\alpha_2}\dots\partial_{x_d}^{\alpha_d}$ and $|\alpha|:=\sum_j \alpha_j$.}
The space $\Cc^k_{b,\Gc}$ is naturally endowed with the $\Cc^k_b$-topology associated to the norm 
$$\|f\|_{\Cc^k}=\sum_{\alpha\in\NN^d,|\alpha|\leq k} \| \partial^\alpha f\|_{L^\infty(X,X)}$$ 
and we recall that this defines a Banach space. 
Then, we endow the whole space $\Cc^k_\Gc$, which includes unbounded functions, with the extended topology. 
\begin{defi}
The {\bf extended topology} of $\Cc^k_\Gc$ is the topology generated by the family of neighborhoods 
\begin{equation}\label{def_N_extended}
\Nc(f_*,\varepsilon):=\{f_*+h,~h\in\Cc^k_{b,\Gc}\text{ and }\|h\|_{\Cc^k}<\varepsilon\}
\end{equation}
where $f_*\in\Cc^k_\Gc$ and $\varepsilon>0$.
\end{defi}
In other words, two (possibly unbounded) functions $f$ and $g$ are close if their difference is bounded and small in the $\Cc^k-$norm. In particular, the extended topology is locally a Banach space and the sequentially closed subsets are closed subsets. Let us also recall the following consequence.
\begin{prop}\label{prop_extended}
For any $k\geq 1$, the space $\Cc^k_\Gc$ is a Baire space: a countable intersection of dense open subsets is dense.   
\end{prop}
\begin{demo}
Let $(\Oc_n)$ be a family of dense open subsets. Let $f_*\in\Cc^k_\Gc$ and let $\Nc$ be a neighborhood of $f_*$. By definition of the topology, there exists $\varepsilon>0$ such that the neighborhood $\Nc(f_*,\varepsilon)$ defined by \eqref{def_N_extended} is included in $\Nc$. Each set $\Oc_n\cap \Nc(f_*,\varepsilon)$ is an open dense subset of the topological space $\Nc(f_*,\varepsilon)$. But this last set is an open subset of a Banach space, so it is a Baire space and in particular $\cap_n \Oc_n \cap \Nc(f_*,\varepsilon)$ is not empty. So $\cap_n\Oc_n\cap \Nc\neq \emptyset$, proving that $\cap_n\Oc_n$ is dense.  
\end{demo}

The above proposition is a guarantee that generic sets provide a relevant notion of ``large'' sets and of ``almost always satisfied'' properties. See Section \ref{sect_discuss_large} for a discussion on other possible notions.
\begin{defi}\label{defi_generic}
A set $\Gg\subset \Cc^k_\Gc$ is {\bf generic} in $\Cc^k_\Gc$ if it contains a countable intersection of dense open sets of $\Cc^k_\Gc$. A property is said to be generic in $\Cc^k_\Gc$ if it is satisfied for a generic set of vector fields.
\end{defi}

\subsection{Coloring and synchrony}\label{section_def_synchrony}
A \emph{coloring} $\bowtie$ of the network is a partition of the set of cells defining a relation of equivalence $c\bowtie c'$ if and only if $c$ and $c'$ are in the same set (i.e. have the same color). 
\begin{defi}
A coloring is {\bf balanced} if, for any cells $c$ and $c'$ with $c\bowtie c'$, there exists an input isomorphism $\beta:I(c)\rightarrow I(c')$ which preserves the colors in the sense that, for all $a\in I(c)$, $T(a)\bowtie T(\beta(a))$.
\end{defi}
Note that two cells having the same color must be input isomorphic and thus of the same cell-type. Moreover, a balanced coloring is such that, if two cells have the same color, then they have the same numbers of input arrows of each type, with the tails of the same color. So the partition defined by the coloring is finer than one defined by the input isomorphisms, which is finer than the ones defined by the cell-types. Also note that a balanced coloring does not necessarily correspond to a symmetry of the network, i.e. to a permutation of the cells under which $\Gc$ and the types are invariant, see \cite{GNS,GS_vs,GS2} or the examples below.
{
One of the main concern in the coupled cell networks is how their cells synchronize, see \cite{Wang} for several examples.}
\begin{defi}
The {\bf synchrony space} of a coloring $\bowtie$ is the subset of $X=\RR^d$ defined by 
$$\Delta_\Join:=\{x\in\RR^d~,~~c\bowtie c'~\Longleftrightarrow ~ x_c=x_{c'}\}.$$
\end{defi}
It is shown in \cite[Theorem 4.3]{GST} and \cite[Theorem 6.5]{SGP} that a coloring is balanced if and only if its synchrony space is invariant for the ODE \eqref{ode} for any admissible map $f\in\Cc^k_\Gc$. In the literature, the equivalence $\Longleftrightarrow$ of the above definition is sometimes an implication $\Rightarrow$. Our present choice is more accurate but note that $\Delta_\Join$ is not a vector space.  
\begin{defi}
Let $x\in\RR^d$ be a state, we define the {\bf synchrony pattern} $\bowtie_{x}$ by 
$$ c\bowtie_x c'~~\Longleftrightarrow ~~x_c=x_{c'}.$$
Let $J\subset\RR$ be an interval of times and let $x(\cdot)\in\Cc^0(J,X)$ be a curve. We define the synchrony pattern $\bowtie_{x,J}$ of the curve, or simply $\bowtie_J$ if there is no confusion, by 
$$ c\bowtie_{x,J} c'~~\Longleftrightarrow ~~x_c(t)=x_{c'}(t)\text{ for all }t\in J.$$
If $J=\RR$, we speak of the ``global'' synchrony pattern and, otherwise, we call such a coloring a ``local'' synchrony pattern. 
In the degenerated case $J=\{t_*\}$, we simply write $\bowtie_{t_*}$ for $\bowtie_{\{t_*\}}$, which is also $\bowtie_{x(t_*)}$. 
\end{defi}

Let us make a small break to consider again the example of Figure \ref{fig_exemple3}. The represented network admits few balanced colorings. 
\begin{itemize}
\item There is of course the trivial one: each cell has a different color. 
\item A more interesting one is the coloring with $c_3\bowtie c_7$ and no other symmetry. Indeed, the cells $c_3$ and $c_7$ have the same type and the same type of input. The condition for $c_3\bowtie c_7$ being balanced is that their inputs have the same color and it holds since $c_3$ and $c_7$ have actually the same input cell $c_5$. It is also clear that if a solution $x(\cdot)$ of \eqref{ode} is such that $x_3(t_*)=x_7(t_*)$ at some time $t_*$, then  $x_3(t)=x_7(t)$ for all $t$ because exchanging $x_3$ and $x_7$ is then harmless for the dynamics of the other cells. 
\item The cells $c_3$ and $c_4$ are input isomorphic. If we start by assuming that $c_3\bowtie c_4$ in a balanced coloring, then, considering their inputs, we must have $c_5\bowtie c_6$ and then, as an input of $c_5$, $c_7$ must be of the same color as one of the inputs cells of $c_6$, that are $c_3$ and $c_4$. Thus, we must have $c_3\bowtie c_4\bowtie c_7$ and $c_5\bowtie c_6$. This yields another balanced coloring. Again, we can see that the corresponding synchrony pattern is natural: there are solutions $x(\cdot)$ of \eqref{ode} such that $x_2(t)=x_3(t)=x_7(t)$ and $x_5(t)=x_6(t)$ because this part of the network is autonomous and the vector field is compatible with these symmetries. Note that this balanced coloring is not completely associated to a symmetry of the graph itself since the cells $c_5$ and $c_6$ are not equivalent in the graph since there is no arrow from $c_6$ to $c_3$.
\item Another balanced coloring that is not a direct symmetry of the graph is the synchrony $c_1\bowtie c_2$. Even if these cells are input isomorphic, they have not the same status in the graph, because even if $a_1$ and $a_2$ are related by the input isomorphism, only $a_1$ is a self-connection. However, once we consider solutions such that $x_1(t)=x_2(t)$, this difference does not matter and both cells follow the same dynamics: this synchrony is preserved. We refer to \cite{GNS,GS_vs,GS2} for more relevant examples of this phenomenon.
\item For a last example, let us consider the coloring $c_2\bowtie c_{9}$. At first sight, it could be possible since the inputs of the cells are of the same type. But, to be balanced, we should have $c_4\bowtie c_7$ and $c_1\bowtie c_{10}$. The last condition is incompatible with a balanced coloring since $c_1$ and $c_{10}$ are not input isomorphic due to the type of the arrow $a_{15}$ which is different from the one of $a_{3}$. Thus, what can be seen as a small defect of the connection $a_{15}$ precludes any synchrony between the left and right part of the network of Figure \ref{fig_exemple3}. 
\end{itemize}

Let us next recall some properties of the synchrony patterns. 
We use the following notion that is classical for relations of equivalence. 
\begin{defi}\label{defi_coarser}
Let $\bowtie$ and $\equiv$ be two coloring. We say that $\bowtie$ is {\bf finer} than $\equiv$ if $c \bowtie c'$ implies $c \equiv c'$. 
We equivalently say that $\equiv$ is {\bf coarser} than $\bowtie$. We say that $\bowtie$ is {\bf strictly finer} than $\equiv$ if $\bowtie$ is finer but not equal to $\equiv$ or equivalently that $\bowtie$ is finer and not coarser than $\equiv$. In the same way, $\equiv$ is strictly coarser than $\bowtie$ if it is coarser but not finer.
\end{defi}
We recall that the synchrony pattern is semi-continuous as already noted in \cite{Stewart}, which can be stated in the following way.
\begin{prop}\label{prop_coarser1}
Let $J\subset\RR$ be an interval of times, let $x(\cdot)\in\Cc^0(J,X)$ be a curve and let $\bowtie_t$ be its synchrony patterns for $t\in J$. For any fixed time $t_*\subset J$, the set 
\begin{equation}\label{defi_J_t}
J_{t_*}:=\{t\in J~,~~\bowtie_t\text{ is finer than }\bowtie_{t_*}\}
\end{equation}
is an open subset of $J$.
\end{prop}
\begin{demo}
Saying that $\bowtie_t$ is finer than $\bowtie_{t_*}$ (i.e. $t$ belongs to $J_{t_*}$) exactly means that $x_c(t)=x_{c'}(t)$ implies $x_c(t_*)=x_{c'}(t_*)$, or equivalently that $x_c(t)\neq x_{c'}(t)$ for each couple of cells where $x_c(t_*)\neq x_{c'}(t_*)$. So the result is a consequence of the continuity of $t\mapsto x(t)$.
\end{demo}

It has also been noted in \cite{Stewart} that $t\mapsto \bowtie_t$ is in general not constant in time: there may exist isolated times where a new symmetry is exceptionally satisfied. We can also switch from one symmetry to another as illustrated in Section 7 of \cite{Stewart}. To work with a constant synchrony pattern, the following result can be useful. It was already mentioned in \cite{Stewart} for example. It is a consequence of the semi-continuity of the synchrony pattern but we provide a basic proof for the sake of completeness.
\begin{prop}\label{prop_coarser2}
Let $J\subset\RR$ be an open non-empty interval of times, let $x(\cdot)\in\Cc^0(J,X)$ be a curve and let $\bowtie_t$ be its synchrony pattern for each $t\in J$. Then, there exists an open non-empty subinterval $J_{0}\subset J$ such that $t\in J_{0}\mapsto~ \bowtie_t$ is constant.
\end{prop}
\begin{demo}
The relation ``is finer than'' is a partial order. Since there is a finite number of possible colorings, there is a time $t_*\in J$ such that there is no $t\in J$ with $\bowtie_t$ strictly finer than $\bowtie_{t_*}$. This is equivalent to saying that, for all $t\in J$, either $\bowtie_{t_*}$ is finer than $\bowtie_t$ or the relations $\bowtie_{t_*}$ and $\bowtie_t$ are not comparable. To construct $t_*$, we can choose a first time $t_0\in J$. If there is no time $t_1$ such that $\bowtie_{t_1}$ is strictly finer than $\bowtie_{t_0}$, then we are done. Otherwise, take such a $t_1$ and continue the process, which will end after a finite number of iterations and provide a suitable time $t_*$.

By Proposition \ref{prop_coarser1}, there exists an open interval $J_{0}$ containing $t_*$ such that, for all $t\in J_{0}$, $\bowtie_{t}$ is finer than $\bowtie_{t_*}$. But, by construction, $\bowtie_{t}$ cannot be strictly finer and so 
$\bowtie_{t}=\bowtie_{t_*}$ for all $t\in J_{0}$.
\end{demo}
{
\begin{remark}
In the above proposition and throughout this article, we adopt the convention that ``$a\in A\mapsto b\in B$'' stands for a function defined from set $A$ to set $B$, mapping the element $a$ to its image $b$.
\end{remark}}

\subsection{Main result}
Having introduced all the necessary notations, we can state rigorously our main result, which is proved in Section \ref{sect_traj}.
\begin{theorem}\label{th_main}
For any $k\geq 1$, there exists a generic set $\Gg\subset \Cc^k_\Gc$ of admissible vector fields such that, for any $f\in\Gg$ and for any solution $x(\cdot)$ of $\dot x(t)=f(x(t))$ in any open time interval $J$, the synchrony pattern $\bowtie_{x,J}$ is balanced.
\end{theorem}
This result means that, for a generic admissible ODE, the only possible synchrony of a solution $x(\cdot)$ are those dictated by the structural constraints of the network and its types. Balanced synchrony can simply correspond to the invariance of the graph under the permutation of its cells. However, more subtle types of balanced synchrony exist, as already noted above, see \cite{GNS,GS_vs,GS2}. Also note that it is always possible that $\bowtie_{x(t)}$ is not balanced for a specific time $t$, but this cannot persist: $\bowtie_{x,J}$ is balanced as soon as $J$ is open. In this sense, the above result extends \cite[Theorem 4.3]{GST} and \cite[Theorem 6.5]{SGP}: a coloring is balanced if and only if its synchrony space is invariant for the ODE \eqref{ode} for \underline{one} generic map $f\in\Cc^k_\Gc$.
In the literature, this kind of results is often restricted to specific solutions as equilibria or periodic orbits. Finally, note that Theorem \ref{th_main} concerns all the solutions, even the ones that are blowing-up in finite time.


\section{Genericity and transversality tools}\label{sect-trans}

In classical problems involving finite-dimensional manifolds, the proofs of the genericity of a property mainly use Sard's Theorem or theorems of transversality similar to the ones of Thom, see \cite{Abraham-Robbin,Sard,Thom1,Thom2}. We will use here a specific result of this type, which goes back to Henry \cite{Henry}. To separate some technical arguments from the proofs of our main results, we also show a ``black-box'' result which is adapted to our context.

\subsection{Henry's theorem}
To extend the transversality theorems to infinite-dimensional context, 
Smale  showed in \cite{Smale} that Sard's Theorem can be extended to Banach
spaces by using the notion of Fredholm operators. Later, Quinn in \cite{Quinn} noted that the notion of left-Fredholm operators is often sufficient. 
See \cite{Henry}, \cite[Section 4.5]{Kato} or \cite{Schechter} for basics properties of Fredholm and semi-Fredholm operators.
\begin{defi}\label{defi_Fred}
Let $X$ and $Y$ be two Banach spaces. A bounded linear operator $L:X\rightarrow
Y$ is a \emph{left-Fredholm operator} if:
\begin{enumerate}[(i)]
 \item its kernel $\Ker(L)$ splits in $X$, i.e. there exists a space $X_1$ such
that $X=X_1\oplus \Ker(L)$ and $X_1$ and $\Ker(L)$ are closed subspaces,
\item its image $\R(L)$ splits in $Y$, i.e. there exists a space $Y_2$ such
that $Y=\R(L)\oplus Y_2$ and both subspaces are closed,
\item its kernel $\Ker(L)$ is finite-dimensional.
\end{enumerate}

If moreover the supplementary space $Y_2$ is also finite-dimensional, then $L$
is called a \emph {Fredholm operator}. The index of $L$ is defined by $\Ind(L)=\dim(\Ker(L))-\dim(Y_2)$ (which is equal to
$-\infty$ if $L$ is not a Fredholm operator).
\end{defi}

\begin{figure}[ht]
\begin{center}
\resizebox{13cm}{!}{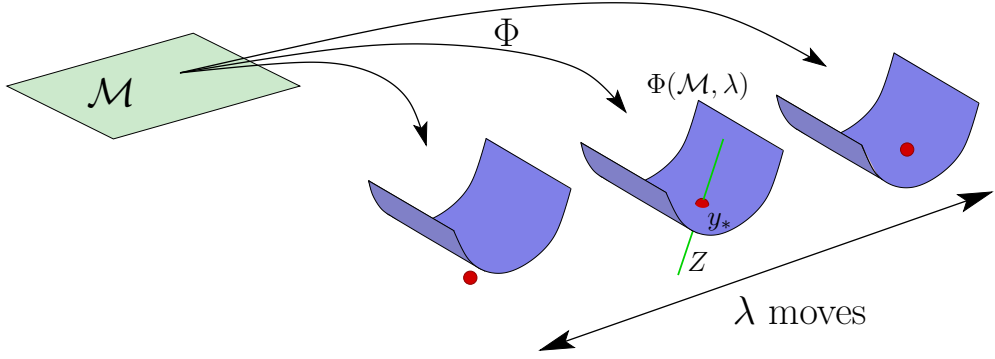} 
\end{center}
\caption{\it A simple illustration of Henry's Theorem. For each $\lambda\in\Lambda$, the function $\Phi(\cdot,\lambda)$ maps a two-dimensional manifold $\Mc$ into a two-dimensional submanifold of $\Nc=\RR^3$. The kernel of $D_x\Phi(x,\lambda)$ is $\{0\}$ and its index is $-1$ because its image is of codimension $1$. Even if $\Phi(\Mc,\lambda)$ may contain a given point $y_*\in\RR^3$ for some specific $\lambda$, if the image of $D\Phi$ contains a direction $Z$ not included in the tangent space $T_{y_*}\Phi(\Mc,\lambda)=\R(D_x\Phi(x,\lambda))$, then $y_*\not\in\Phi(\Mc,\lambda)$ for a generic $\lambda$.}\label{fig_Henry} 
\end{figure}

Following Smale arguments, one can extend the classical transversality theorems to Banach manifolds. There exist many different
versions of this kind of theorems in Banach manifolds, often called Sard-Smale
theorems, see for example \cite{Abraham-Robbin}, \cite{Henry} or \cite{Saut-Temam}. In the present paper, we use the following version, proved by Henry (see Figure \ref{fig_Henry} for illustration). It corresponds to Theorem 5.4 of \cite{Henry} with Assumption 2.$\beta$ and the use of the remark following the statement for Assumption 3. We also refer to another version of the theorem and its proof in \cite{RJ-prevalence}. 
\begin{theorem}[Henry's Theorem]\label{th_Henry}
$~$\\
Let $\Mc$, $\Lambda$ and $\Nc$ be three Banach manifolds. Let $\Phi:\Mc\times
\Lambda \longrightarrow \Nc$ be a map of class $\Cc^1$ and $y_*$ be a point of
$\Nc$.\\ 
We assume that :
\begin{description}
\item{(i)} $\forall (x,\lambda)\in \Phi^{-1}(y_*)$, 
$D_x\Phi(x,\lambda):T_x\Mc \rightarrow T_{y_*}(\Nc)$ is a left-Fredholm 
  operator with negative index,
\item{(ii)} $\forall (x,\lambda)\in \Phi^{-1}(y_*)$,  the image of the total
derivative $D\Phi(x,\lambda):T_x\Nc\times
T_\lambda\Lambda\rightarrow T_{y_*}\Nc$ contains a finite-dimensional subspace $Z$
such that $Z\cap \R(D_x\Phi(x,\lambda))=\{0\}$ and the dimension of $Z$ is
strictly larger than the one of $\Ker(D_x\Phi(x,\lambda))$,
\item{(iii)} $\Mc\times\Lambda$ is separable.
\end{description}
Then there exists a generic subset $\Gg$ of $\Lambda$ such
that, for any $\lambda_0 \in\Gg$, $y_*$ is not in the image of the map $x\mapsto
\Phi(x,\lambda_0)$. 
\end{theorem}

In our applications of Theorem \ref{th_Henry}, the operator $D_x\Phi$ can be
split as $D_x\Phi=L+K$, where $K$ is a compact operator and $L$ is a simple
operator, for which Hypotheses (i) and (ii) are more easily checked. The following classical result (see for example \cite[Section 4.5]{Kato} or \cite[Theorem 5.22]{Schechter}) implies that it is indeed sufficient to check (i) for the operator $L$.
\begin{prop}\label{prop_Fred_1}
Let $X$ and $Y$ be two Banach spaces. Let $L:X\rightarrow Y$ be a
left-Fredholm  operator and let $K\in\Lc(X,Y)$ be a compact operator. Then,
$L+K$ is a  left-Fredholm map with the same index as the one
of $L$. 
\end{prop} 
{
To check Hypothesis (ii) of Theorem \ref{th_Henry} for an operator $D_x\Phi$ of the form $L+K$, we will use the following property. It requires the construction of a suitable space $Z$ related only to $L$ and this is often simpler than constructing a space for the entire operator $L+K$. Note that the following proposition is more general than the one used in \cite{RJ}, since $Z$ is not included in $Y_2$, but roughly said in a sector around $Y_2$.}
\begin{prop}\label{prop_Fred_3}
Let $L:X\rightarrow Y$ be a left-Fredholm map and $Y_2$ be as in
Definition \ref{defi_Fred}. Let $\pr$ be the continuous projection on $Y_2$ corresponding to the splitting $Y=\R(L)\oplus Y_2$. Assume that $Z$ is a subspace of  $Y$ such that there exists $\kappa>0$ such that 
$$\forall z\in Z~,~~\|\pr(z)\|\geq \kappa\|z\|.$$
Then, for any compact operator $K:X\rightarrow Y$, the subspace $Z\cap \R(L+K)$ is finite-dimensional.
\end{prop}
\begin{demo}
Let $(z_n)$ be any bounded sequence in $Z\cap \R(L+K)$. Due to Proposition \ref{prop_Fred_1}, we know that $L+K$ is a left-Fredholm operator. In particular, there exists a closed splitting $X=\Ker(L+K)\oplus X_1$. The Banach
isomorphism theorem implies that $(L+K)$ restricted to $X_1$ admits a bounded inverse $T:\R(L+K)\rightarrow X_1$. We set $x_n=T(z_n)$, defining a bounded sequence of $X$. Applying the projection $\pr$, we obtain that $\pr (z_n)=\pr(L+K)x_n=\pr K x_n$. Since $(x_n)$ is bounded and $K$ is compact, this shows that we can extract a converging subsequence $(\pr (z_{\varphi(n)}))$. Since $z_n$ belongs to $Z$, we have that 
$$\forall n,m\in\NN~,~~\|\pr (z_{\varphi(n)})-\pr (z_{\varphi(m)})\|\geq \kappa\|z_{\varphi(n)}- z_{\varphi(m)}\|.$$
This implies that $(z_{\varphi(n)})$ is also a Cauchy sequence, which is thus convergent. To summarize, we have shown that any bounded sequence of $Z\cap \R(L+K)$ admits a convergent subsequence. This implies that $Z\cap \R(L+K)$ is a finite-dimensional space.
\end{demo}

\subsection{An adapted black-box}

Theorem \ref{th_Henry} will be used several times during our proofs. To check Hypotheses (i) and (ii), certain key arguments would be repeated. Some of them are purely a matter of functional analysis and topology. To avoid making the future proofs heavy, we decide to gather most of the topological arguments in the present section. To this end, we are introducing a general framework that will suit our applications. The resulting Proposition \ref{prop_bb} below is referred as a ``black-box'' in the sense that, during the future proofs, we will apply it directly without concern for the functional analysis techniques contained within its proof.

\medskip

We use the notations of Section \ref{section_notations}.
Let $(\sigma,\tau)\in\RR^2$ with $\sigma<\tau$ be two times. Let $c$ be a cell, we denote $\tilde X_c=\{x\in X,x_{c'}=0\text{ if }c'\neq c\}$ the canonical embedding of the state space $X_c$ in $X$. Let $\hat \pr:X \rightarrow \tilde X_c$ be a surjective projection. { To be able to consider two slightly different frameworks without repeating two propositions and two very similar proofs, we 
introduce the subspace $\Pc$ of $\Cc^1([\sigma,\tau],\RR^d)$ as either
$$\Pc:=\{x(\cdot)\in \Cc^1([\sigma,\tau],\RR^d)~,~~\hat \pr(x(t))=0\text{ for all }t\in[\sigma,\tau]\}$$
or 
$$\Pc:=\{x(\cdot)\in \Cc^1([\sigma,\tau],\RR^d)~,~~t\in[\sigma,\tau] \mapsto \hat \pr(x(t))\text{ is constant}\}.$$
In the following, $\Pc$ is one the above sets with $\hat \pr$ a suitable chosen projection.}
For example, if we have two cells $c$ and $c'$ with the same state space, choosing $\hat \pr(x)=(0,\ldots,0,x_c-x_{c'},0,\ldots)$ and the first case generates the subspace $\Pc$ of curves having the states in cells $c$ and $c'$ are equal for all times (notice that we do not assume $\hat \pr$ to be necessarily the canonical coordinate projection). 
To give another example, choosing $\hat \pr(x)=(0,\ldots,0,x_c,0,\ldots)$ and the second case generates the subspace $\Pc$ of curves for which the state in cell $c$ is stationary. 
\begin{prop}\label{prop_bb}
Consider the above framework and let $\Oc$ be an open subset of $\Pc$. Let $k\geq 1$.
Assume that, for any admissible vector field $f\in\Cc^k_\Gc$ and any solution $t
\in[\sigma,\tau]\mapsto x(t)$ of the ODE 
\begin{equation}\label{ode_bb}
\dot x(t)=f(x(t))~~~~t\in[\sigma,\tau]
\end{equation}
with $x(\cdot)\in\Oc$, there exists a space $G\subset \Cc^k_\Gc$ of admissible vector fields such that 
\begin{enumerate}[(a)]
 \item the space $Z:=\{g\circ x,~g\in G\}$ is a subspace of $\Cc^0([\sigma,\tau],X)$ with infinite dimension,  
 \item there exists $\kappa >0$ such that, for all $z\in Z$, $\|\hat\pr z\|_{L^\infty([\sigma,\tau],\tilde X_c)}\geq \kappa \|z\|_{L^\infty([\sigma,\tau],X)}$.
\end{enumerate}

Then, there exists a generic set $\Gg\subset\Cc^k_\Gc$ of admissible vector fields such that, for all $f\in\Gg$, there is no solution $x(\cdot)$ of $t\in[\sigma,\tau]\mapsto x(t)$ of the ODE \eqref{ode_bb} belonging to the set $\Oc$.
\end{prop}
\begin{demo} 
We aim at describing the generic set $\Gg$ as the countable intersection of dense open subsets. We first introduce these suitable subsets $\Gg_{n,m}$ as follows.
By construction, $\Pc$ is a closed subspace of $\Cc^1([\sigma,\tau],X)$. Moreover, it admits a closed complementary space $\Qc$ in $\Cc^1([\sigma,\tau],X)$. Indeed, since $\hat\pr$ is a projector onto $\tilde X_c$, $\Ker (\hat\pr)\oplus \R(\hat\pr)=\Ker (\hat\pr)\oplus \tilde X_c=X$ and we can take 
$$\Qc:=\{x(\cdot)\in \Cc^1([\sigma,\tau],X)~,~~x(t)\in \tilde X_c\text{ for all }t\in[\sigma,\tau]\}$$
in the first case or  
$$\Qc:=\{x(\cdot)\in \Cc^1([\sigma,\tau],X)~,~~x(t)\in \tilde X_c\text{ for all }t\in[\sigma,\tau]\text{ and }x(\sigma)=0\}$$
in the second case.
We also recall that the open subset $\Oc\subset\Pc$ can be written as a countable union $\Oc=\cup_{n\in\NN} \Fc_n$ of closed sets $\Fc_n$ of $\Pc$. Indeed, we can choose 
$$\Fc_n:=\big\{f\in \Pc, d(f,\Pc\setminus\Oc)\geq \frac 1n\big\}$$
where $d(f,\Pc\setminus\Oc)=\inf\{ \|f-g\|_{\Cc^k} , {g\in \Pc\text{ and }g\not\in\Oc}\}$ (actually, the fact that open sets are countable union of closed sets is a general property of metric spaces). Finally note that the sets $\Fc_n$ are also closed in the whole space $\Cc^1([\sigma,\tau],X)$ since $\Pc$ is closed. 
{
For all integers $n$ and $m$, we set 
\begin{align*}
\Gg_{n,m} := \{ &f\in\Cc^k_\Gc\text{ such that there is no solution }t\in[\sigma,\tau]\mapsto x(t)\\ &\text{of the ODE \eqref{ode_bb} with }\|x(\cdot)\|_{\Cc^1}\leq m\text{ and }x(\cdot)\in\Fc_n\}. 
\end{align*}
It remains to prove that these sets are open and dense. Proving the density is the more complicated task and will require Theorem \ref{th_Henry}. We split the proof in several steps.}

\medskip

{\noindent\it$\rhd$ Step 1:  The set $\Gg_{n,m}$ is open in $\Cc^k_\Gc$.}\\ 
Consider a function $f$ in $\Gg_{n,m}$. By definition of the extended topology, to prove the openness of $\Gg_{n,m}$, it is sufficient to show that, for all $g\in\Cc^k_{b,\Gc}$ small enough, $f+g$ belongs to $\Gc_{n,m}$. We argue by contradiction: assume that there exists a sequence $(g_p)$ of functions converging to $0$ in $\Cc^k_{b,\Gc}$ such that $f+g_p\not\in\Gg_{n,m}$. By definition, for all $p$, there exists a solution $t\in[\sigma,\tau]\mapsto x_p(t)$ of the ODE 
\begin{equation}\label{ode_bb2}
\dot x_p(t)=f(x_p(t))+g_p(x_p(t))~~~~t\in[\sigma,\tau]
\end{equation}
satisfying $\|x_p(\cdot)\|_{\Cc^1}\leq m$ and $x_p(\cdot)\in\Fc_n$. We know that the sequence $(x_p)$ is bounded in $\Cc^1$. Using \eqref{ode_bb2}, the classical bootstrap argument shows that it is also bounded in $\Cc^2([\sigma,\tau],\RR^d)$. Using Ascoli's theorem, we have that $(x_p)$ is compact in $\Cc^1$ and, up to renumbering the sequence, we can assume that $(x_p)$ converge to a function $x_\infty$ in $\Cc^1([\sigma,\tau],\RR^d)$. Since $\Fc_n$ is closed in $\Cc^1([\sigma,\tau],\RR^d)$, the function $x_\infty$ also belongs to $\Fc_n$. Also note that the uniform bound $\|x_p(\cdot)\|_{\Cc^1}\leq m$ yields $\|x_\infty(\cdot)\|_{\Cc^1}\leq m$. Finally, passing to the limit in \eqref{ode_bb2}, we get that $x_\infty(t)$ is a solution of the ODE $\dot x_\infty(t)=f(x_\infty(t))$. This limit would contradict the fact that $f$ belongs to $\Gg_{n,m}$. We conclude that there exists a small neighborhood of $f$ included in $\Gg_{n,m}$.

\medskip

{\noindent\it$\rhd$ Step 2:  The set $\Gg_{n,m}$ is dense in $\Cc^k_\Gc$.}\\
Let $f_*$ be any function of $\Cc^k_\Gc$ and $\Vc$ be a neighborhood of $f_*$ in $\Cc^k_\Gc$. Let $\Lambda$ be the set of admissible vector fields $h\in\Cc^k_\Gc$ with support included in the ball of radius $m+2$. We aim at finding $h\in\Lambda$ as small as wanted such that $f=f_*+h$ belongs to $\Gg_{n,m}$. By definition of the extended topology, choosing a small enough $h$ will ensure that $f$ belongs to $\Vc$ and this will prove the density of $\Gg_{n,m}$. Note that the perturbations $h$ have compact support, which is technically important because it ensures that $\Lambda$ is a separable Banach space, which is not the case of the whole space $\Cc^k_\Gc$. Since we only consider in $\Gg_{n,m}$ solutions $x(\cdot)$ with values in the ball of radius $m$, this restriction on the support will be harmless. Next, we set $\Nc=\Cc^0([\sigma,\tau],\RR^d)$ and $y_*$ be the zero function of $\Nc$. Finally, we recall that $\Oc$, introduced in the statement of Proposition \ref{prop_bb}, is a submanifold of $\Cc^1([\sigma,\tau],\RR^d)$ since it is an open subset of a closed subspace with closed complementary space. We set $\Mc:=\Oc\cap B(0,m+1)$ the set of functions $x(\cdot)$ in $\Oc$ with $\sup_{t\in[\sigma,\tau]} |x(t)|+|\dot x(t)|<m+1$. We introduce the function $\Phi\in\Cc^1(\Mc\times\Lambda,\Nc)$ defined by
$$\Phi(x,h)=\frac{\d x}{\d t}(\cdot)-f_*(x(\cdot))-h(x(\cdot))~.$$
Note that $\Phi(x,h)=0$ exactly means that $x$ is a solution of $\dot x(\cdot)=(f_*+h)(x(\cdot))$ belonging to the ball $B(0,m+1)$ and to the set $\Oc$. We aim at applying Theorem \ref{th_Henry} to the above framework. We have that $\Mc$ and $\Lambda$ are separable Banach manifolds. Assume for the moment that Hypotheses (i) and (ii) of Theorem \ref{th_Henry} hold. Then, for a generic vector field $h \in \Lambda$, the function $0$ is not in the image of $\Phi(\cdot,h)$. This means that, for a generic $h$, there is no solution of $\dot x(\cdot)=(f_*+h)(x(\cdot))$ belonging to the ball $B(0,m+1)$ and to the set $\Oc$. This implies in particular that, for a generic $h$, $f_*+h$ belongs to $\Gg_{n,m}$. Since $h$ can be taken in a generic set, it can be chosen as small as needed, proving the density of $\Gg_{n,m}$ in a neighborhood of $f_*$.

\medskip

{\noindent\it$\rhd$ Step 3: $D_x\Phi$ is a left-Fredholm map.}\\
Let us check that Hypothesis (i) of Theorem \ref{th_Henry} holds, with the notations introduced in Step 2. Remember that $(x,h)\in\Phi^{-1}(y_*)$ exactly means that $x\in\Mc$ is a solution of the ODE $\dot x(\cdot)=(f_*+h)(x(\cdot)):=f(x(\cdot))$. 
Since $\Oc$ is an open subset of the vector space $\Pc$, all the tangent spaces of $\Oc$ are equal to $\Pc$. The function $\Phi$ is of class $\Cc^1$ and we have
$$D\Phi(x,h).(\xi,g)=\frac{\d \xi}{\d t}(\cdot)
-D(f_*+h)(x(\cdot)).\xi(\cdot)- g(x(\cdot))~.$$
The map $L:\xi \in \Pc \mapsto \frac{\d \xi}{\d t} \in \Cc^0([\sigma,\tau],X)$ is a left-Fredholm function. Indeed, its kernel is the set of the constant functions of $\Pc$ and therefore is of dimension at most $d$ and admits a closed complementary set consisting of the functions of $\Pc$ vanishing at $\sigma$. 
In both options of definition of $\Pc$, the range of $L$ is the set 
\begin{equation}\label{eq_demo_bb2}
R(L)=\{y(\cdot)\in \Cc^0([\sigma,\tau],X)~,~~\hat \pr(y(t))=0\text{ for all }t\in[\sigma,\tau]\}
\end{equation}
which admits a closed complementary space being
\begin{equation}\label{eq_demo_bb}
\{y(\cdot)\in \Cc^0([\sigma,\tau],X)~,~~y(t)\in \tilde X_c\text{ for all }t\in[\sigma,\tau]\}
\end{equation}
Ascoli's theorem yields that the embedding of $\Cc^1([\sigma,\tau],\RR^d)$ in $\Cc^0([\sigma,\tau],\RR^d)$ is compact. Since $X=\RR^d$ and since $D(f_*+h)(x(t))$ belongs to $\Cc^0([\sigma,\tau],\Lc(\RR^d))$, the map $$K:\xi \in \Cc^1([\sigma,\tau],X)\mapsto D(f_*+h)(x(t)).\xi \in\Cc^0([\sigma,\tau],X)$$ is compact. Applying Proposition \ref{prop_Fred_1}, we obtain that
$D_x\Phi$ is a left-Fredholm map. 

\medskip

{\noindent\it$\rhd$ Step 4: there is enough freedom to construct suitable perturbations.}\\
Again, with the notations introduced in Step 2, $(x,h)\in\Phi^{-1}(y_*)$ exactly means that $x\in\Mc$ is a solution of the ODE $\dot x(\cdot)=(f_*+h)(x(\cdot))$ and, thus, we can use hypotheses (a) and (b) of Proposition \ref{prop_bb} to check that hypotheses (ii) of Theorem \ref{th_Henry} holds. 
We have just seen that $D_x\Phi$ can be written $L+K$ as in Proposition \ref{prop_Fred_3} with $Y_2$ given by \eqref{eq_demo_bb}. 
In particular, we note that the projection on $Y_2$ along $R(L)$ is the projection $\pr$ given by
$$\pr~:~y(\cdot)\in\Cc^0([\sigma,\tau],X)~\longmapsto (t\mapsto \hat\pr(y(t))).$$
Note that $g\circ x=-D\Phi(x,h).(0,g)$, meaning that the infinite dimensional space $Z$ provided by the hypotheses of Proposition \ref{prop_bb} is a subspace of the image of the total derivative $D\Phi(x,h)$. 
We can apply Proposition \ref{prop_Fred_3} to this setting: the part of $Z$ belonging to the range of $L+K=D_x\Phi$ is finite-dimensional. Thus, even if we get rid off this part, we can still find a subspace $\tilde Z$ of $Z$, with large enough dimension, such that Hypothesis (ii) of Theorem \ref{th_Henry} holds. 
To conclude, note that, even if the previous arguments seem to show rigorously that we can apply Theorem \ref{th_Henry}, there is still a small technical gap: the fields $g$ provided by Assumption (ii) are not necessarily with compact supports as the above definition of the parameter space $\Lambda$ requires. However, the considered solution $x$ belongs to $\Mc$ and is therefore valued in $B(0,m+1)$. So we can smoothly truncate the fields $g$ to obtain fields belonging to $\Lambda$ without changing the values of $g\circ x$. 

\medskip

{\noindent\it$\rhd$ Step 5: conclusion}\\
Gathering all the previous steps, we have shown that $\Gg_{n,m}$ is an open subset and that we can apply Theorem \ref{th_Henry} to obtain that $\Gg_{n,m}$ is also dense. Thus, $\Gg:=\cap_{n,m\in\NN} \Gg_{n,m}$ is a generic subset of admissible vector fields $\Cc^k_G$ such that, for any $f\in\Gg$, there is no solution of the ODE \eqref{ode_bb} belonging to $\Oc$.
\end{demo}


\section{A strategy for constructing perturbations}\label{section_perturb}

To  apply our ``black box'' Proposition \ref{prop_bb}, we need to be able to construct a family of admissible vector fields $g$ satisfying its assumptions (a) and (b). Recall that $\hat \pr$ is a projection onto the state of a cell $c$. So we need to be able to construct with an infinite-dimensional freedom (Assumption (a)) perturbations $g$ valued in a cone oriented along the state space $X_c$ (Assumption (b)). Consider a solution $t\mapsto x(t)$ of the ODE $\dot x(t)=f(x(t))$ on a time interval $J\subset\RR$. Let us focus on the cell $c$ and its input cells $T(I(c))$ and write more shortly $T:=T(I(c))$. Assume that $(x_c(\cdot),x_T(\cdot))$ is not stationary and, up to a restriction of $J$, assume that $\frac{\d~}{\d t}(x_c(\cdot),x_T(\cdot))$ never vanish on $J$. Then, up to make $J$ even smaller, we have that $t\in J\mapsto (x_c(\cdot),x_T(\cdot))$ is diffeomorphic to a curve and we can construct functions $g_c$ such that $t\mapsto g_c(x_c(t),x_T(t))$ covers any $\Cc^1-$curve in $X_c$. This is basically how the constructions of \cite{RJ} are made, since this article considers networks without constraint of symmetry. In the case of networks with types and their associated ODE, we can use several ideas, mostly coming from \cite{GST,Stewart,Stewart_eq} or other previous works.

\medskip 

\begin{figure}[ht]
\begin{center}
\resizebox{13cm}{!}{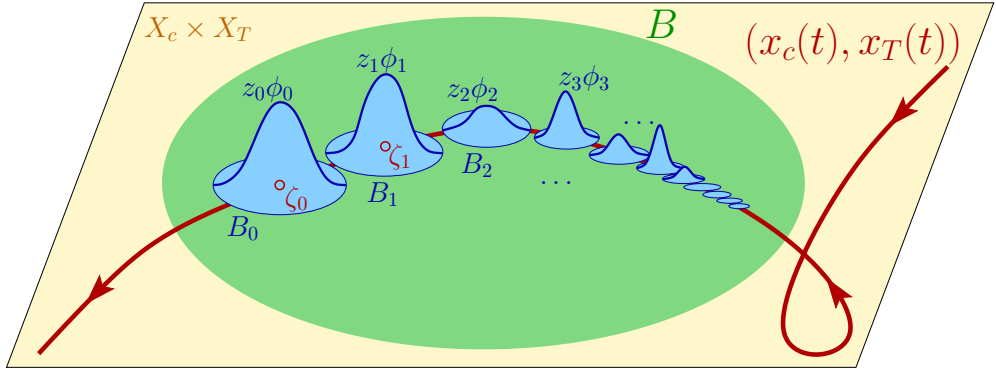} 
\end{center}
\caption{\it To construct a suitable perturbation, we focus on a ball $B$ where $t\mapsto (x_c(t),x_T(t))$ is a bijective curve. Then a function $\phi$ is generated by the combination of bump functions $\phi_n$ with disjoint supports. The choice of the amplitudes $z_n$ of the bumps provides an infinite-dimensional freedom. Moreover, the resulting function $\phi=\sum z_n\phi$ reaches its maximum along the curve $t\mapsto (x_c(t),x_T(t))$.}\label{fig_perturbation} 
\end{figure}

\noindent{$\rhd$ \it Trick 1: Generating a suitable infinite dimensional space.}\\
We start with  the assumption that $t\in J\mapsto (x_c(t),x_T(t))$ is a $\Cc^1-$diffeomorphism describing a $\Cc^1-$curve included in a given ball $B$ of $X_c\times X_T$. We choose a sequence $t_n=t_*+\eta 2^{-n}$ ($n\in\NN$) with $t_*\in J$ and with $\eta>0$ small enough such that $t_*+\eta$ also belongs to $J$, implying that each $t_n$ belongs to $J$. Each point $\zeta_n:=(x_c(t_n),x_T(t_n))$ is at positive distance of the other ones and there exists radius $\varepsilon_n$ such that each ball $B_n:=B(\zeta_n,\varepsilon_n)\subset X_c\times X_T$ is disjoint from the others and included in $B$. Define $\phi_n\in\Cc^k(X_c\times X_T,\RR_+)$ as a smooth non-negative ``bump'' function with support in $B_n$, reaching its maximal value at $\zeta_n$ and normalized by $\|\phi_n\|_{\Cc^k}=1$ (see Figure \ref{fig_perturbation} below for an illustration). For any sequence $(z_n)\in\ell^1(\NN)$, we set $\phi=\sum_n z_n \phi_n$, which is a well-defined function of $\Cc^k(X_c\times X_T,\RR_+)$. Note that each $t\in J\mapsto \phi_n(x_c(t),x_T(t))$ has a support disjoint from the others and is not zero. This yields two important properties for us.
First, the space of functions $t\in J\mapsto \phi_n(x_c(t),x_T(t))$ generated by all the choices of $(z_n)\in\ell^1(\NN)$ is an infinite-dimensional vector space. Second, the maximal value of $\phi$ is exactly the maximum of $n\mapsto \|z_n\phi_n\|_\infty$ and also the maximum of the real function $t \in J\mapsto \phi(x_c(t),x_T(t))$.

\medskip 

\noindent{$\rhd$ \it Trick 2: Symmetrization.}\\
The symmetrization is a trick already used in \cite{Maxime,Stewart_eq,Stewart}. Consider a function $\phi$ generated by the previous process and fix a unit vector $y\in X_c$. Then, using the previous notations and the ones of Section \ref{section_notations}, we set
\begin{equation}\label{eq_section_perturb1}
g_c(x_c,x_T):=\sum_{\beta\in B(c,c)} \phi(x_c,\beta^* x_T)y.
\end{equation}
The obvious but important remarks are that $g_c$ only depends on the input cells of $c$ and $g_c$ is symmetric with respect all the admissible permutations of the inputs of the cell $c$. Then, for any cell $c'$ which is input equivalent to $c$, we set 
\begin{equation}\label{eq_section_perturb2}
g_{c'}(\zeta,\xi):=  g_{c}(\zeta,\beta^* \xi)~~~\text{ for some }\beta\in B(c,c').
\end{equation}
This definition is independent of the choice of $\beta$ since, if $\beta'$ is another input isomorphism, then $\beta^{-1}\circ \beta'$ belongs to $B(c,c)$ and $g_c$ is invariant under actions of $B(c,c)$. For other cells $c''$, we set $g_{c''}\equiv 0$. Following the above remarks, we can check that this construction generates an admissible vector field $g$. 

\medskip 

\noindent{$\rhd$ \it Trick 3: Avoiding destructive interactions.}\\
During the above symmetrization process, we would like to avoid that the different terms of the sums could balance out. For example, we have an infinite-dimensional freedom to construct the function $\phi$ of the first step. But the summation could reduce this freedom: e.g. there exist infinite-dimensional spaces of functions of $\Cc^0(\RR,\RR)$ such that the space of their images by the symmetrization $f\mapsto f(\cdot)+f(-\cdot)$ is finite-dimensional. Another important point concerns the estimation of the maximum of the function $t \in J\mapsto \|g_c(x_c(t),x_T(t))\|$ appearing in Assumption (b) of Proposition \ref{prop_bb}. If we know that the maximum of the real function $t \in J\mapsto \phi(x_c(t),x_T(t))$ is exactly $\|\phi\|_{L^\infty}$, destructive interactions in the symmetrization process may reduce this maximum. 

To avoid this kind of problems, we show the following useful lemma.
\begin{lemma}\label{lemme1}
Let $J\subset\RR$ be an open non-empty interval of times and let $c$ be a cell, having $T:=T(I(c))$ as input cells. Let $x(\cdot)\in\Cc^1(J,\RR^d)$ be a curve such that $\frac{\d~}{\d t}(x_c(t),x_T(t))$ never vanishes on $J$.
Then, there is an open non-empty subinterval $J_{*}\subset J$ and an open ball $B\in X_c\times X_T$ such that:
\begin{enumerate}[(i)]
\item the curve $t\in J_{*} \mapsto (x_c(t),x_T(t))$ is a $\Cc^1-$diffeomorphim on its image, which is included in $B$,
\item for all input isomorphisms $\beta\in B(c,c)$ of the cell $c$: 
\begin{itemize}
\item either $(x_c(t),\beta^*x_T(t))=(x_c(t),x_T(t))$ for all $t\in J_{*}$, 
\item or $(x_c(t),\beta^*x_T(t))\not\in B$ for any $t\in J_{*}$.
\end{itemize}
\end{enumerate}
\end{lemma}
\begin{demo}
We use Proposition \ref{prop_coarser2}: there exists an open interval $J_0\subset J$ such that the synchrony pattern of $x(t)$ is constant. Thus, for any given input isomorphism $\beta\in B(c,c)$, if $\beta^* x_T(t)=x_T(t)$ for some $t\in J_0$ then  $\beta^* x_T(t)=x_T(t)$ for all $t\in J_0$. 
Choose a time $t_0\in J_0$, $\eta>0$ small enough such that $(t_0-\eta,t_0+\eta)\subset J_0$. Since, by assumption, $\frac{\d~}{\d t}(x_c(t_0),x_T(t_0))\neq 0$, up to choose $\eta>0$ smaller, the curve $t\in (t_0-\eta,t_0+\eta) \mapsto (x_c(t),x_T(t))$ is a $\Cc^1-$diffeomorphim on its image. By construction, for any $\beta\in B(c,c)$, either $\beta^* x_T(t)=x_T(t)$ for all $t\in (t_0-\eta,t_0+\eta)$ or this equality never holds. In this last case, we have, in particular, $\beta^* x_T(t_0)\neq x_T(t_0)$ and we can make $\eta>0$ smaller to find $\varepsilon>0$ such that the ball $B((x_c(t_0),x_T(t_0)),\varepsilon)$ contains the image of $t\in (t_0-\eta,t_0+\eta) \mapsto (x_c(t),x_T(t))$ but no image of $t\in (t_0-\eta,t_0+\eta) \mapsto (x_c(t),\beta^* x_T(t))$. There is only a finite number of input isomorphisms $\beta\in B(c,c)$, so we can repeat the restriction process for all of those such that $\beta^* x_T(t)= x_T(t)$ never holds. In the end, we obtain an open interval $ J_*\subset J$ and a small ball $B=B((x_c(t_0),x_T(t_0)),\varepsilon)$ such that the statement holds.   
\end{demo}

The conclusion of Lemma \ref{lemme1} is illustrated in Figure \ref{fig_lemme1} below. Its interest is the following. If $\phi$ is supported in a small enough ball $B$, the terms of the sum \eqref{eq_section_perturb1} are either zero or the same for all $t\in J_{*}$, where $J_*$ is provided by Lemma \ref{lemme1}. Then $t\in J_{*} \mapsto g_c(x_c(t),x_T(t))$ is a simple multiple of $\phi(x_c(\cdot),x_T(\cdot))$. This type of arguments is not new, see for example Lemma 7.2 of \cite{Stewart_eq}, which concerns the case of equilibria.

\begin{figure}[ht]
\begin{center}
\resizebox{11cm}{!}{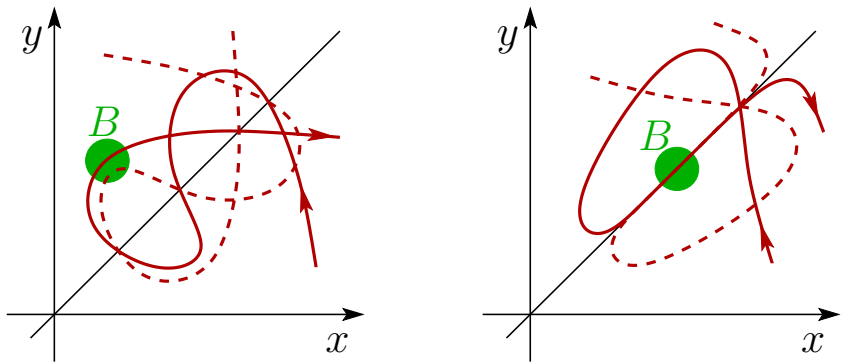} 
\end{center}
\caption{\it An simplified illustration of Lemma \ref{lemme1} for a curve $t\mapsto (x(t),y(t))$ (plain line) and its symmetry by the permutation $\beta^*(x,y)=(y,x)$ (dashed line). In the situation on the left, the curve satisfies the symmetry only at exceptional times and we can find a ball $B$ satisfying the second point of (ii) of Lemma \ref{lemme1}. Note that $B$ may contain points of the curve $t\mapsto \beta^* (x(t),y(t))$, but not in the interval $J_*$. In the situation on the right, the curve satisfies the symmetry for an open interval. Then, it is possible to choose a ball $B$ satisfying the first point of (ii) of Lemma \ref{lemme1}.}\label{fig_lemme1} 
\end{figure}


\section{Generic synchrony of equilibrium points}\label{sect_eq}

The strategy explained in Section \ref{section_perturb} deals with solutions for which each cell has non-stationary inputs. This means that the case of stationary solutions has to be studied separately. This is the purpose of the present section.

We use the notations of Section \ref{section_def_synchrony}. Our first result proves that the synchrony pattern of equilibrium points is almost always balanced. It generalizes  \cite[Theorem 7.6]{GST} or \cite{Stewart_eq} since no hyperbolicity is required in our result. Note that the genericity of hyperbolicity is not guaranteed because of the presence of strong symmetries. So skipping this assumption is meaningful. 

\begin{prop}[Generic balanced synchrony for equilibria]\label{prop_synchrony_eq}
Let $k\geq 1$. 
There exists a generic set $\Gg\subset \Cc^k_\Gc$ of admissible vector fields such that, for all $f\in\Gg$ and for all $x\in\RR^d$, if $f(x)=0$ then the synchrony pattern $\bowtie_x$ of $x$ is balanced. 
\end{prop}
\begin{demo}
The core of the proof is an application of Theorem \ref{th_Henry}. Due to the stationary inputs, we cannot use the black-box Proposition \ref{prop_bb} directly but we will follow the same ideas.

We fix a coloring $\bowtie$ of the cells and $m>0$. Assume that $\bowtie$ is not balanced and let $\Oc_m\subset \Cc^k_\Gc$ be the set of admissible vector fields such that there is no $x\in X$ with $f(x)=0$, $\|x\|\leq m$ and 
\begin{equation}\label{eq_demo_equi}
\left\{\begin{array}{ll} 
x_c=x_{c'} & \text{ if }c\bowtie c'\\
\|x_c-x_{c'}\|\geq 1/m & \text{ if }X_c=X_{c'}\text{ without having }c\bowtie c'.
\end{array}\right. 
\end{equation}

First, we easily check that $\Oc_m$ is open. Indeed, assume that there is a sequence $(f_n)\subset \Cc^k_\Gc\setminus \Oc_{m}$ converging to $f_\infty$ and a bounded sequence $(x_n)$ such that $f_{n}(x_n)=0$, $\|x_n\|\leq m$ and \eqref{eq_demo_equi} holds. We can extract a subsequence $(x_{\varphi(n)})$ that converges to some $x_\infty$. Passing to the limit, we can check that $f_\infty(x_\infty)=0$ with $\|x_\infty\|\leq m$ and \eqref{eq_demo_equi} also holds for the limit $x_\infty$. This shows that the complement of $\Oc_{m}$ is closed. 

To prove the density, let us fix $f\in \Cc^k_\Gc$. We consider perturbations of $f$ of the form $f+g$ with $g$ supported in $B(0,m+2)$. The interest of such perturbations is that $$\Lambda:=\Cc^k_\Gc\cap \Cc^k(\overline B(0,m+2),\RR^d)$$ is a separable Banach space, contrarily to $\Cc^k_\Gc$ endowed with its extended topology.  We introduce a new equivalence between cells denoted $\equiv$ by setting $c\equiv c'$ if and only if $c\bowtie c'$ and there exists an input isomorphism $\beta:I(c)\rightarrow I(c')$ preserving the coloring $\bowtie$, that is that $T(a)\bowtie T(\beta(a))$ for any input arrow $a$ of $c$. We denote by $C$ a set containing exactly one unique element of each class of equivalence of $\equiv$. We can explain the motivation of $\equiv$ as follows. Let $x$ be a point with synchrony pattern $\bowtie_x=\bowtie$ and let $g\in\Cc^k_\Gc$. If $c$ and $c'$ are two input isomorphic cells, then, up to the pullback by an isomorphism $\beta$, $g_c$ and $g_{c'}$ are the same. So if we know that $c \equiv {c'}$, the information $g_{c'}(x)=0$ is redundant with $g_c(x)=0$ since the value of $x$ in the cells $c$ and $c'$ and their input cells are equal. So $C$ is exactly a minimal set of cells such that $g_C(x)=0$ implies $g(x)=0$ for all $x$ with $\bowtie_x=\bowtie$. 
We apply Henry's Theorem with 
$$\Mc=\{x\in\RR^d,~\bowtie_x=\bowtie\text{ and }\|x\|<m+1\}\text{ and }\Nc=X_C.$$
We set $y_*=0_\Nc$ and 
$$\Phi:(x,g)\in\Mc\times \Lambda ~\longmapsto~ (f+g)_C(x)\in X_C.$$
Assume for the moment that the hypotheses of Theorem \ref{th_Henry} hold in this framework. Then, its conclusion yields the existence of a generic set of perturbations $g$ such that $f+g$ 
has no solution of $(f+g)_C(x)=0$ with $\bowtie_x=\bowtie$ and $\|x\|<m+1$. By construction of $C$, this means that there is no zero of $f+g$ with $\bowtie_x=\bowtie$ and $\|x\|<m+1$ and this implies in particular that $f+g$ belongs to $\Oc_m$. Since $g$ is generic in $\Lambda$, it can be taken as small as wanted and this shows that $\Oc_m$ is dense.

It remains to check the hypotheses of Theorem \ref{th_Henry}. Note that $\Mc$ and $\Nc$ are finite-dimensional and that Hypothesis (iii) is obvious. Due to the rank–nullity theorem, the index of $D_x\Phi$ is $\dim (\Mc)-\dim (\Nc)$. The dimension of $\Mc$ is $d_K$ where $K$ is a set of representative cells of the equivalence relation $\bowtie$, whereas the dimension of $\Nc$ is $d_C$ where $C$ is a set of representative cells of the equivalence relation $\equiv$. By definition, $c\equiv c'$ implies $c\bowtie c'$ so $\equiv$ is finer than $\bowtie$ and $\dim (\Mc)\leq \dim (\Nc)$. To obtain the strict equality, we recall that $\bowtie$ is assumed to be unbalanced: there exists a least a couple of cells with $c\bowtie c'$ but without any input isomorphism $\beta\in B(c,c')$ preserving $\bowtie$. The cells $c$ and $c'$ have the same color for $\bowtie$ but are not in the same class of equivalence for $\equiv$. So $\equiv$ is strictly finer than $\bowtie$ and the index of $D_x\Phi$, equal to $\dim (\Mc)-\dim (\Nc)$, is negative. To check Assumption (ii) of Theorem \ref{th_Henry}, it is sufficient to show that, for any $x\in\Mc$ and $g\in\Lambda$, $D\Phi(x,g)$ is surjective onto $X_C=T_{y_*}\Nc$ because of $\dim (\Mc)<\dim (\Nc)$. The derivative $D_\lambda\Phi(x,g).\tilde g$ is equal to $\tilde g_C(x)$ and it is sufficient to show that we can choose $\tilde g_C(x)$ freely. To this end, we  construct it on each cell $c_1$, $c_2$, \ldots, $c_p$ of $C$ step by step. Assume that $\tilde  g_C(x)$ is constructed in $c_1$, \ldots, $c_j$ and consider $c_{j+1}$. There are three cases:
\begin{enumerate}[1)]
\item there is no cell $c_i$ with $i\leq j$ such that $c_i$ and $c_{j+1}$ are input isomorphic. Then, the symmetries of the network require no constraint on $\tilde g_{c_{j+1}}$ and $\tilde g_{c_{j+1}}(x)$ can be freely chosen.
\item there exists $i\leq j$ such that $c_i$ and $c_{j+1}$ are input isomorphic, $c_i\bowtie c_{j+1}$ and there exists an input isomorphism between $c_i$ and $c_{j+1}$ preserving $\bowtie$. This exactly means that $c_i \equiv c_{j+1}$ and these cells cannot both belong to $C$. So this case is actually impossible by construction.
\item there exists $i\leq j$ such that $c_i$ and $c_{j+1}$ are input isomorphic but  either $c_i$ and $c_{j+1}$ have not the same color for $\bowtie$, or there is no input isomorphism preserving $\bowtie$. In both case, this means that, for all input isomorphism $\beta\in B(c_{j+1},c_i)$, we always have $(x_{c_i}, x_{T(I(c_i))})\neq (x_{c_{j+1}},\beta^* x_{T(I(c_{j+1}))})$. So, even if $g_{c_{j+1}}$ is constrained by the symmetries of the network, we can modify it in a small neighborhood of $(x_{c_{j+1}},x_{T(I(c_{j+1}))})$ and symmetrize it as in Trick 2 of Section \ref{section_perturb}, without changing its values on any $(x_{c_i},x_{T(I(c_i))})$ with $i\leq j$, that is in the already considered cells.
Thus, we can choose the value of $\tilde g_c(x_{c_{j+1}},x_{T(I(c_{j+1}))})$ independently from the previously constructed values.
\end{enumerate}

To conclude, the previous arguments show that we can apply Theorem \ref{th_Henry}
to show that $\Oc_m$ is dense. Thus $\Gg_{\bowtie}:=\cap_{m\in\NN\setminus\{0\}} \Oc_m$ is a generic set consisting of functions $f$ having no zero with synchrony pattern $\bowtie$. 
There is only a finite number of choices for the unbalanced coloring $\bowtie$. So intersecting the corresponding sets $\Gg_{\bowtie}$ provides the desired generic set $\Gg$ of the statement. 
\end{demo}

Note that one of the core arguments of the above proof is the fact that $\dim (\Mc)<\dim (\Nc)$. This can be associated to the method of ``overdetermined constraints'' of Stewart, see \cite{Stewart_eq,Stewart}.  

As already said, we must deal with the case of stationary inputs separately. The previous result concerns equilibrium points where the state is stationary in all cells. But there may exist solutions that are stationary in some cells but not all. To avoid the degenerated situations, we prove the result below. It can be considered as the observation of stationary cells since it ensures that, generically, a state in a cell is constant if and only if all its inputs are constant, and so are the inputs of its inputs etc. Note that it is still possible to have a stationary input to a non-stationary cell and that this situation has no reason to be exceptional.
\begin{prop}[Observation of stationary cells]\label{prop_rigid_constant}
Let $k\geq 1$. There exists a generic set $\Gg\subset \Cc^k_\Gc$ of admissible vector fields such that the following property hold. For any $f\in\Gg$, for any solution $x(\cdot)$ of $\dot x(t)=f(x(t))$ in any open time interval $J$ and for any cell $c$, if $t\in J \mapsto x_c(t)$ is constant, then in any input cell $c'\in T(I(c))$ of $c$, $t\in J \mapsto x_{c'}(t)$ is also constant.
\end{prop}
\begin{demo}
First consider two times $\sigma<\tau$. Let us fix a given cell $c$ and set $T:=T(I(c))$ being the set of input cells of $c$. We apply our ``black-box'' Proposition \ref{prop_bb} with the projection $\hat \pr(x)=x_c$, the subspace 
$$\Pc=\{x(\cdot)\in\Cc^1([\sigma,\tau],\RR^d)~,~~t\mapsto \hat \pr(x(t)):=x_c(t)\text{ is constant}\}$$
and its open subset
$$\Oc=\{x(\cdot)\in\Pc~,~~t\mapsto x_{T}(t)\text{ is not constant}\}.$$
Let us fix any solution $x(\cdot)$ of $\dot x(t)=f(x(t))$ belonging to $\Oc$. 
We need to construct the subspace $G$ as required by the hypotheses (a) and (b) of Proposition \ref{prop_bb}. To this end, we use the tricks introduced in Section \ref{section_perturb}. 

Since $x_{T}$ is not constant and of class $\Cc^1$, there is a interval of times $J$ where $t\in J\mapsto \frac{\d~}{\d t}x_{T}(t)$ never vanishes and this is a fortiori the same for $\frac{\d~}{\d t}(x_c(t),x_{T}(t))$. We apply Lemma \ref{lemme1}: there exists $J_{*}\subset J$ and a ball $B\in X_c\times X_T$ such that the properties (i) and (ii) of Lemma \ref{lemme1} hold. We introduce the space $G$ as the space of admissible vector fields generated by the two first tricks of Section \ref{section_perturb}. More precisely:
\begin{enumerate}[1)]
\item We construct a family $(\phi_n)\subset\Cc^k(B,\RR_+)$ of bump functions with disjoint supports whose maximum is reached at a point of the curve $t\in J_{*}\mapsto (x_c(t),x_{T}(t))$ as in Trick 1 of Section \ref{section_perturb}. This family generates an infinite-dimensional space of functions $\phi$ supported in $B$, with $\max |\phi|$ reached somewhere along the curve $t\in J_{*}\mapsto (x_c(t),x_{T}(t))$. Moreover, the family of functions $t\in J_{*} \mapsto \phi(x_c(t),x_{T}(t))$ is also infinite-dimensional.
\item From each function $\phi$ of the above family, we construct an admissible vector field $g\in\Cc^k_\Gc$ by following the symmetrization process of Trick 2 of Section \ref{section_perturb}. Denote by $G$ the space of all the functions $g$ obtained by this construction. 
\end{enumerate}
It remains to check that $G$ satisfies Assumption (a) and (b) of Proposition \ref{prop_bb}. A key remark is that, due to Lemma \ref{lemme1}, a field $g\in G$ is such that, in the cell $c$, 
$$\forall t\in J_{*}~,~~g_c(x(t)):=\sum_{\beta\in B(c,c)} \phi(x_c(t),\beta^* x_{T}(t))y=K \phi(x_c(t),x_{T}(t))y,$$
where $K\geq 1$ is the number of input isomorphisms such that the first option of (ii) of Lemma \ref{lemme1} holds (which is at least $1$ because of $\beta=\id$). As a first consequence, the space $Z:=\{g\circ x,g\in G\}$ is infinite-dimensional. Indeed, the space of functions obtained as $t\in J_{*}\mapsto \phi(x_c(t),x_{T}(t))$ is already infinite-dimensional and this property remains true when extending in other cells and extending the time interval $J_{*}$ to $J$. A second consequence is that $\sup_{t\in J_{*}} \|g_c(x(t))\|=K\|\phi\|_{L^\infty}$ because the functions $\phi$ reached their maxima along the curve $t\in J_{*}\mapsto (x_c(t),x_{T}(t))$. So, remembering that $\hat \pr$ is here the canonical projection on the component of the cell $c$, we can check that Assumption (b) of Proposition \ref{prop_bb} holds. Indeed, for each $z=g\circ x$ with $g\in G$, we have 
$$\|\hat \pr z\|_{L^\infty}=\sup_{t\in J} \|g_c(x(t))\|\geq \sup_{t\in J_{*}} \|g_c(x(t))\| = K\|\phi\|_{L^\infty}$$
and, by construction of $g$ following Trick 2 of Section \ref{section_perturb}, 
$$\|z\|_{L^\infty}\leq \|g\|_{L^\infty} = \max_{\text{cells }d} \|g_d\|_{L^\infty} = \|g_c\|_{L^\infty} \leq \sharp B(c,c) \|\phi\|_{L^\infty} .$$
As a conclusion, we get 
$$ \|z\|_{L^\infty}\leq \frac{\sharp B(c,c)}{K} \|\hat \pr z\|_{L^\infty}$$
and we can apply our ``black-box'' Proposition \ref{prop_bb}: there exists a generic set $\Gg_{c,\sigma,\tau}\subset \Cc^k_\Gc$ of admissible vector fields such that there is no solution of the ODE $\dot x(t)=f(x(t))$ belonging to $\Oc$, that is such that $x_c(\cdot)$ is constant in $[\sigma,\tau]$ but $x_T(\cdot)$ is not constant. In other words, for all $f\in \Gg_{c,\sigma,\tau}$, if a solution is such that $x_c(\cdot)$ is constant in $[\sigma,\tau]$, $x_{c'}(\cdot)$ is also constant in $[\sigma,\tau]$ for all the input cells $c'$ of $c$. 

To finish the proof of Proposition \ref{prop_rigid_constant}, it remains to set 
$$\Gg:=\bigcap_{\text{cells }c} \left( \bigcap_{ (\sigma,\tau)\in\QQ^2 \text{ with }\sigma<\tau} \Gg_{c,\sigma,\tau}\right).$$
The set $\Gg$ is generic in $\Cc^k_\Gc$ as a countable intersection of generic sets (remember that $\Cc^k_\Gc$ is a Baire space as noted in Section \ref{section_def_ODE}). 
\end{demo}


\section{Generic synchrony of trajectories}\label{sect_traj}
In the previous section, we have been focusing on stationary cells. Here, we consider cells with a non-stationary state and we conclude by proving our main result Theorem \ref{th_main}.

\begin{prop}[Balanced synchrony of non-constant cells]\label{prop_traj}
Let $k\geq 1$. There exists a generic set $\Gg\subset \Cc^k_\Gc$ of admissible vector fields such that the following property hold for any $f\in\Gg$. For any solution $x(\cdot)$ of $\dot x(t)=f(x(t))$ in any open time interval $J$ and for any cells $c\neq c'$, if $x_c(t)=x_{c'}(t)$ and $\dot x_c(t)\neq 0$ for all $t\in J$, then $c$ and $c'$ are input isomorphic and there exists an input isomorphism $\beta\in B(c,c')$ such that the input cells satisfy $x_{T(I(c'))}(t)=\beta^* x_{T(I(c))}(t)$ for at least one $t\in J$.
\end{prop}
\begin{demo}
The strategy of the proof is very similar to the one of Proposition \ref{prop_rigid_constant}. First consider two times $\sigma<\tau$. Fix two cells $c\neq c'$ and set $T:=T(I(c))$ and $T':=T(I(c'))$ the set of their input cells. We apply our ``black-box'' Proposition \ref{prop_bb} with the projection $\hat \pr(x)=x_c-x_{c'}$, the subspace 
\begin{align*}
\Pc&:=\{x(\cdot)\in\Cc^1([\sigma,\tau],\RR^d)~,~~\hat \pr(x(\cdot))\equiv 0\}\\
&=\{x(\cdot)\in\Cc^1([\sigma,\tau],\RR^d)~,~~\forall t\in [\sigma,\tau],~x_c(t)=x_{c'}(t)\}
\end{align*}
and its open subset
$$\Oc=\{x(\cdot)\in\Pc~,~~\forall t\in [\sigma,\tau],~\dot x_c(t)\neq 0\text{ and }\forall \beta\in B(c',c),~\beta^* x_{T'}(t)\neq x_{T}(t)\}.$$
Remember that $B(c',c)$ is the (possibly empty) set of input isomorphisms from $c'$ to $c$, that are in finite number. Using in addition the compactness of $[\sigma,\tau]$, it is easy to show that $\Oc$ is indeed an open subset of $\Pc$. 

Let us fix a solution $x(\cdot)$ of $\dot x(t)=f(x(t))$ belonging to $\Oc$. We only have to explain how to construct the subspace $G$ as required by the hypotheses (a) and (b) of Proposition \ref{prop_bb}. As in the proof of Proposition \ref{prop_rigid_constant}, we use the tricks introduced in Section \ref{section_perturb}. First, we use Lemma \ref{lemme1} to obtain a subinterval $J_{*}\subset (\sigma,\tau)$ and a ball $B$ as in (i) and (ii) of Lemma \ref{lemme1}.
Now, we consider all the possible images of the inputs of the other cell $c'$ by input isomorphisms.
For each $\beta\in B(c',c)$, by definition of $\Oc$, we have $\beta^*x_{T'}(t)\neq x_{T}(t)$ for all $t\in J_*$. Thus, up to choose $B$ and $J_*$ smaller, we can assume that $t\in  J_*\mapsto (x_{c'}(t),\beta^* x_{T'}(t))$ is always outside $B$. 

Let $G$ be the space of admissible vector fields generated by the two first tricks of Section \ref{section_perturb}. More precisely, copying the proof of Proposition \ref{prop_rigid_constant}:
\begin{enumerate}[1)]
\item We construct a family $(\phi_n)\subset\Cc^k(B,\RR_+)$ of bump functions with disjoint supports whose maximum is reached at a point of the curve $t\in  J_* \mapsto (x_c(t),x_{T}(t))$ as in Trick 1 of Section \ref{section_perturb}. It generates an infinite-dimensional space of functions $\phi$ supported in $B$, with $\max |\phi|$ reached somewhere along the curve $t\in  J_* \mapsto (x_c(t),x_{T}(t))$. We again note that the family of functions $t\in  J_* \mapsto \phi(x_c(t),x_{T}(t))$ is also infinite-dimensional.
\item From each function $\phi$ of the above family, we construct an admissible vector field $g\in\Cc^k_\Gc$ by following the symmetrization process of Trick 2 of Section \ref{section_perturb}: a unit vector $y\in X_c$ is chosen and $g$ is defined by \eqref{eq_section_perturb1} and \eqref{eq_section_perturb2}. Denote by $G$ the space of all the functions $g$ obtained by this construction. 
\end{enumerate}
Checking that $G$ satisfies Assumption (a) of Proposition \ref{prop_bb} follows the same argument as in the proof of Proposition \ref{prop_rigid_constant}: due to Lemma \ref{lemme1}, the symmetrization process does not destroy the fact that the constructed family of functions $t\in  J_* \mapsto g_c(x_c(t),x_{T}(t))$ is infinite-dimensional because
$$\forall t\in  J_*~,~~g_c(x(t))=K\phi(x_c(t),x_T(t))y$$
for some positive integer $K$. So the family $t\in [\sigma,\tau] \mapsto (g\circ x)(t)$ is also infinite-dimensional. To check Assumption (b), the important remark is that the curve $t\in  J_*\mapsto (x_{c'}(t), x_{T'}(t))$ is outside the support of any $g\in G$ because its image by any input isomorphism is outside $B$. Thus, for any $t\in  J_*$ and $g\in G$, we have 
$$\hat \pr (g\circ x)(t)=(g_c\circ x)(t)-(g_{c'}\circ x)(t)=(g_c\circ x)(t)=K\phi(x_c(t),x_T(t))y.$$
Thus, for all $z=g\circ x$ with $g\in G$, we have $\|\hat \pr z\|_{L^\infty} \geq K\|\phi\|_{L^\infty}$. On the other hand,  
$$\|z\|_{L^\infty}\leq \|g\|_{L^\infty} = \max_{\text{cells }d} \|g_d\|_{L^\infty} = \|g_c\|_{L^\infty} \leq \sharp B(c,c) \|\phi\|_{L^\infty}.$$
So we can apply our ``black-box'' Proposition \ref{prop_bb}: there exists a generic set $\Gg_{c,c',\sigma,\tau}\subset \Cc^k_\Gc$ of admissible vector fields such that there is no solution of the ODE $\dot x(t)=f(x(t))$ belonging to $\Oc$. This exactly means that the conclusion of Proposition \ref{prop_traj} hold for the particular choice of cells and for the time interval $[\sigma,\tau]$. To finish the proof of the proposition, it remains to intersect all the above sets for all the possible couple of cells and all the times $\sigma<\tau$ with $(\sigma,\tau)\in\QQ^2$. Indeed, we note that if the conclusion of Proposition \ref{prop_traj} hold for a time interval $[\sigma,\tau]$, then it holds for any time interval $J$ containing $[\sigma,\tau]$. 
\end{demo}

We are now able to prove our main result: for a generic vector field, the synchrony patterns of the solutions are always balanced.

\medskip

{ \noindent \emph{\textbf{Proof of Theorem \ref{th_main}:}}}
We construct the suitable generic set by intersecting all the ones provided by the previous propositions. More precisely:
\begin{itemize}
\item We denote by $\Gg_1$ the generic set of Proposition \ref{prop_rigid_constant}, implying that, for all $f\in\Gg_1$, constant cells of a solution $x(t)$ of $\dot x(t)=f(x(t))$ must have constant inputs.
\item We denote by $\Gg_2$ the generic set of Proposition \ref{prop_traj}, implying that, for all $f\in\Gg_2$, synchronous non-constant cells of a solution $x(t)$ of $\dot x(t)=f(x(t))$ must be input equivalent with, at least, a punctual synchrony of the inputs. 
\item We use Proposition \ref{prop_synchrony_eq} in a more subtle way. For each set $C$ of cells, we consider all its indirect inputs, that are the cells from which we can follow a sequence of arrows to arrive at a cell $c\in C$. We consider the subgraph $\Gc_C$ constructed by restricting the whole graph $\Gc$ to these cells and the arrows linking them: the cells of $\Gc_C$ are the cells of $C$ and all their indirect inputs cells in $\Gc$, the arrows of $\Gc_C$ are all the arrows of $\Gc$ whose head is a cell kept in $\Gc_C$ and the types of the elements of $\Gc_C$ are the types of the corresponding elements in $\Gc$. By removing the necessary cells and arrows to obtain $\Gc_C$ from $\Gc$, we transform any admissible $f\in\Cc^k_\Gc$ to a map $\tilde f\in\Cc^k_{\Gc_C}$, which is also admissible. Indeed, all the inputs of any cell $c'\in\Gc_C$ are included in $\Gc_C$ by construction. We can apply Proposition \ref{prop_synchrony_eq} to all these subgraphs $\Gc_C$ obtained from all the possible sets of cells $C$: this provides generic sets $\Gg_C\subset\Cc^k_{\Gc}$ of admissible vector fields such that, for any $f\in\Gc_C$, the synchrony pattern of a solution of $\dot x(t)=f(x(t))$ being constant in the independent subgraph $\Gc_C$ must have a balanced synchrony pattern in this subgraph.
\item Then, we define the generic set 
$$\Gg:=\Gg_1\cap\Gg_2\cap \left(\bigcap_{\text{subset of cells }C }\Gg_C \right).$$
\end{itemize}
It remains to check that $\Gg$ is as claimed by the statement of Theorem \ref{th_main}. Let $f\in\Gg$, let $J$ be any open time interval and let $x(\cdot)$ be a solution of $\dot x(t)=f(x(t))$ for $t\in J$. We consider the synchrony pattern $\bowtie_{x,J}$ and we have to show that it is balanced. 

We first use Proposition \ref{prop_coarser2} to find a subinterval $J_{0}\subset J$ such that $t\in J_{0} \mapsto \bowtie_{x(t)}$ is constant, simply denote this synchrony pattern by $\bowtie$. If we consider a cell $c$ where $t\in J_{0}\mapsto x_c(t)$ is not constant, then, up to restricting $J_{0}$, we can assume that $\dot x_c(\cdot)$ never vanishes in $J_{0}$. Doing this possible restriction successively in all the cells, we can assume that $J_{0}$ is small enough such that, in any cell $c$, either $x_c(\cdot)$ is constant or $\dot x_c(\cdot)$ never vanishes.

Let $c$ and $c'$ be two cells such that $c\bowtie c'$, meaning that $x_c(t)=x_{c'}(t)$ for all $t\in J_{0}$. There are two cases:
\begin{enumerate}[(i)]
\item either $t\in J_{0}\mapsto x_c(\cdot)$ is constant (and so is $x_{c'}(\cdot)$). In this case, Proposition \ref{prop_rigid_constant} recursively implies that all the indirect inputs of the cells $c$ and $c'$ are constant. Then, using Proposition \ref{prop_synchrony_eq} in the corresponding subgraph $\Gc_C$ with $C=\{c,c'\}$, we obtain the existence of an input isomorphism $\beta\in B(c,c')$ such that $\beta^* x_{I(c)}=x_{I(c')}$.
\item or $t\in J_{0}\mapsto \dot x_c(\cdot)$ never vanishes. In this case, we can apply Proposition \ref{prop_traj}: there exists an input isomorphism $\beta\in B(c,c')$ such that $\beta^* x_{I(c)}(t)=x_{I(c')}(t)$ for at least one $t\in J_{0}$. But, by construction, $t\in J_{0} \mapsto \bowtie_{x(t)}$ is constant and so $\beta^* x_{I(c)}(t)=x_{I(c')}(t)$ for all $t\in J_{0}$.
\end{enumerate}
We have just proven that $\bowtie_{x,J_{0}}$ is balanced. It remains to recall that \cite[Theorem 4.3]{GST} yields that the synchrony space of a balanced coloring is invariant for the ODE $\dot x=f(x)$ if $f\in\Cc^k_\Gc$. So the synchrony pattern of $x$ being balanced in the time interval $J_{0}$ remains the same for all time and we have in particular $\bowtie_{x,J}=\bowtie_{x,J_{0}}$, proving that the synchrony pattern $\bowtie_{x,J}$ is balanced.
{\hfill$\square$\\}


\section{Further results and discussions}
\subsection{Rigid patterns of synchrony}\label{sect_discuss_rigid}
Several articles as \cite{GRW,GRW2,GS3,GST,Stewart,SGP,SP} study the rigidity of the synchrony patterns. Let $f\in\Cc^k_\Gc$ be an admissible vector field and $\Nc$ a neighborhood of $0$ in $\Cc^k_\Gc$. Let $J$ a time interval and consider a family of solutions $x_g(t)$ of the ODE $\dot x_g(t)=(f+g)(x_g(t))$ with $t\in J$ and $g\in\Nc$. We assume that $x_g$ depends continuously on $g$ in the sense that $g\in\Nc\mapsto x_g\in\Cc^1(J,X)$ is a continuous map. Classical and important examples of such families of solutions $x_g$ are: families of simple (or even hyperbolic) equilibrium points and families of simple (or even hyperbolic) periodic orbits of $f+g$. 
\begin{defi}\label{defi_rigid}
Consider the above framework. We say that the synchrony pattern $\bowtie_{x_g,J}$, defined in Section \ref{section_def_synchrony}, is {\bf rigid} if, up to choose the neighborhood $\Nc$ smaller, $\bowtie_{x_g,J}=\bowtie_{x_0,J}$ for all perturbations $g\in\Nc$.
\end{defi}
We deduce from our main result the following rigidity property.
\begin{coro}\label{coro_rigid}
Consider a family of solutions $x_g(\cdot)$ as above. If its synchrony pattern is rigid, then $\bowtie_{x_0,J}$ must be balanced. 
\end{coro}
\begin{demo}
This is an obvious consequence of Theorem \ref{th_main}. Indeed, consider that $x_g(\cdot)$ has a rigid synchrony pattern, that is that $\bowtie_{x_g,J}$ is constant for all $f+g$ in a small neighborhood of $f$. By Theorem \ref{th_main}, we can find a perturbation $g$ such that all the synchrony patterns of solutions of the ODE with vector field $f+g$ are balanced. A fortiori, this is the case of $\bowtie_{x_g,J}$ and, by rigidity, of $\bowtie_{x_0,J}$ and all the synchrony patterns of $x_{\tilde g}(\cdot)$ for any $\tilde g\in\Nc$.  
\end{demo}

Applying Corollary \ref{coro_rigid} to a family $g\mapsto x_g$ of hyperbolic equilibrium points, we recover \cite[Theorem 7.6]{GST} or the main result of \cite{Stewart_eq}. 
Applying Corollary \ref{coro_rigid} to a family $g\mapsto x_g(\cdot)$ of (strongly) hyperbolic periodic orbits, we recover Theorem 9.2 of \cite{Stewart}. 
But note that we do not actually need any hyperbolicity for applying Corollary \ref{coro_rigid}. This may be of importance since the Kupka-Smale property is not known to be generic for coupled cell networks with types. It may happens that hyperbolicity fails to be generic in this type of systems.

\subsection{The doubled network and the phase-shift synchrony}
The doubled network is a simple but powerful trick, used in \cite{GRW,GRW2,Stewart_eq,Stewart} for example. Consider a network $\Gc$ with types and its associated space $\Cc^k_\Gc$ of vector fields. Define the doubled network $2\Gc$ as the network consisting in two copies $\Gc_1$ and $\Gc_2$ of $\Gc$, these two copies being disconnected but having exactly the same type of cells and arrows. The important remark is that $F\in\Cc^k_{2\Gc}$ is an admissible vector field for the doubled graph if and only if $F$ is two copies $(f_1,f_2):=(f,f)$ of the same vector field $f\in\Cc^k_\Gc$, where $f_i$ is $F$ restricted to $\Gc_i$. Indeed, the cells and arrows of $\Gc_2$ are copies of the ones of $\Gc_1$ with the same types and \eqref{hyp_f_symmetry} yields that $f_1=f_2$ (with the obvious identifications).  It is also clear that this doubling identification is compatible with the topology endowing the vector fields and thus maps generic sets to generic sets. We refer to \cite[Lemma 4.3]{GRW} or \cite[Section 11]{Stewart}.

We can use the doubled network as in \cite{Stewart} to obtain results on the phase-shift synchrony. Indeed, consider a solution $x(\cdot)$ on the original network $\Gc$ and associate to it the solution $X(\cdot)=(x(\cdot),x(\cdot+\theta))$ on the doubled network $2\Gc$. The synchrony pattern not only identify the equalities $x_c(t)=x_{c'}(t)$ in the original network but also equalities $X_{c_1}(t)=X_{c_2}(t)$ with $c_1$ a cell in $\Gc_1$ and $c_2$ the same cell in $\Gc_2$. This means that the synchrony pattern also detects the phase-shift $x_c(t)=x_c(t+\theta)$ in cells $c$ of the original network. 

This simple trick has the following direct consequence: generically, if two cells have the same dynamics but shifted in time, then they must be input equivalent and in particular they are of the same type.
\begin{coro}[Rigid phase property]\label{coro_phase_shift_1}
Let $\Gc$ be a network with types and let $k\geq 1$. There exists a generic set $\Gg\subset\Cc^k_\Gc$ such that if $f\in\Gg$, if $x(\cdot)$ is a solution of $\dot x(t)=f(x(t))$ in a open time interval $(\sigma,\tau)$ and if there are two cells $c$ and $c'$ and $\theta\in (0,\tau-\sigma)$ such that $x_c(t)=x_{c'}(t+\theta)$ for all $t\in (\sigma,\tau-\theta)$, then $c$ and $c'$ are input equivalent and there exists $\beta\in B(c,c')$ such that $\beta^*x_{T(I(c))}(\cdot )=x_{T(I(c'))}(\cdot+\theta)$.  
\end{coro}
\begin{demo}
We use the doubled network $2\Gc$ and apply Theorem \ref{th_main} to it. The generic set $\Gg$ of Corollary \ref{coro_phase_shift_1} is the restriction of the ones of Theorem \ref{th_main} to the copy $\Gc_1$ of $\Gc$. If $x(\cdot)$ is a solution of $\dot x(t)=f(x(t))$ for $t\in(\sigma,\tau)$, then $X=(x(\cdot),x(\cdot+\theta))$ is a solution of $\dot X(t)=F(X(t))$ for all $t\in (\sigma,\tau-\theta)$ where $F=(f,f)$ is the doubled vector field. Then Corollary \ref{coro_phase_shift_1} simply follows from the fact that the synchrony pattern of $X$ must be balanced, applied to the cells $c_1$ and $c'_2$ of $2\Gc$.  
\end{demo}

With Corollary \ref{coro_phase_shift_1}, we recover Theorem 11.1 and Corollary 11.2 of \cite{Stewart}. { As an interesting application, this kind of shifted dynamics can be seen in animal locomotion, where the phase shifts in the periodic motions of small parts (legs, fins, cilia\dots) generate a global undulation, see \cite{SW}. The previous result suggests that the groups of neurons involved in this kind of locomotion must be symmetric.}

We also have the following important consequence.
\begin{coro}[Observation of oscillations]\label{coro_phase_shift_2}
Let $\Gc$ be a network with types and let $k\geq 1$. There exists a generic set $\Gg\subset\Cc^k_\Gc$ such that the following holds for all $f\in\Gg$ and all global solution $x(\cdot)$ of $\dot x(t)=f(x(t))$, $t\in \RR$. 

If there is a cell $c$ such that $x_c(\cdot)$ is periodic of period $\theta>0$, then for all the input cells $c'$ of $c$, $x_{c'}(\cdot)$ is periodic of period $\theta'\in \theta\NN$. This also hold for the indirect input cells of $c$.  

Moreover, assume that $\Gc$ is transitive, i.e. if each cell is an indirect input of all the others. If in one cell $c$ the state $x_c(\cdot)$ is periodic with minimal period $\theta$, then the whole state $x(\cdot)$ is periodic with minimal period $\theta'=K\theta$ for some integer $K\geq 1$.
\end{coro}
\begin{demo}
We apply Corollary \ref{coro_phase_shift_1} and consider the same generic set $\Gg$. Let $f\in\Cc^k_\Gc$ and $x(\cdot)$ be a global solution of the ODE. Assume that $x_c(\cdot)$ is periodic of period $\theta>0$ in a cell $c$ and denote $T:=T(I(c))$ the input cells of $c$. We consider the conclusion of Corollary \ref{coro_phase_shift_1} with $c=c'$: there exists an input isomorphism $\beta\in B(c,c)$, permuting the inputs of $c$, such that $\beta^*x_T(t)=x_T(t+\theta)$. If $c'\in T$ is invariant by $\beta$, then $x_{c'}(\cdot)$ is $\theta-$periodic. If not, $c'$ belongs to a cycle of the permutation $\beta$ and there is a power $K$ such that $\beta^K(c')=c'$. Since  $(\beta^K)^*x_T(t)=x_T(t+K\theta)$, this implies that $c'$ is $K\theta-$periodic. Obviously, we can iterate the argument to reach all the indirect inputs of $c$. 

If $\Gc$ is transitive, the above argument shows that the state in all the cells is periodic with a period being a multiple of $\theta$. Thus, the whole state $x(\cdot)$ is periodic with period $T$ being the lowest common multiple of all the cells with period $K\theta$. Note that $T$ is not necessarily the minimal period of $x(\cdot)$ but, the state in the cell $c$ having minimal period $\theta$, $x(\cdot)$ must have a minimal period which is anyway a multiple of $\theta$. 
\end{demo}
{ Note that Corollary \ref{coro_phase_shift_2} considers global solutions. The reason is not only to properly define the periodicity. Indeed, consider a solution existing for a finite time interval $[0,T]$, and assume that, in a cell $c$, $x_c(t+\theta)=x_c(t)$ for some $\theta\in (0,T)$ and for all $t$ where it makes sense. In our proof, we show that the input of $c$ has a period $K\theta$ with an integer $K\geq 1$. But if $K\theta$ is greater than $T$, this becomes an empty condition. However, with a careful look to the proof of Corollary \ref{coro_phase_shift_2}, we see that there actually exists an integer $K\geq 1$, which only depends on the geometry of the network, such that if a solution is defined in $[0,T]$ and is $\theta$-periodic in a cell $c$ with $K\theta<T$, then the solution is periodic in all the indirect input cells of $c$.}

\subsection{Other notions of large sets}\label{sect_discuss_large}
In this article, we use the genericity to give a meaning to the notion
of ``almost every vector fields on $X$'', see Definition \ref{defi_generic}. The genericity is a  classical and well-accepted notion of ``large sets'' in infinite-dimensional spaces and it is sufficient to obtain the density of a property and the rigidity of patterns as in Corollary \ref{coro_rigid}. However, we can argue that, in finite dimensional spaces, it is possible to have generic sets with zero measure. This is troublesome since both acceptable notion of ``large sets'' are then contradictory. 

For the interested reader, we recall that there are other notions of ``large sets'' in Banach spaces that are different from the genericity and more related to the Lebesgue measure in finite dimension. If $X$ is a Banach space, \cite{Christ} introduces the notion of Haar-nul set: a Borel set $B$ of $X$ is said \emph{Haar-nul} if there exists a finite non-negative measure $\mu\not\equiv 0$ with compact
support such that for all $x\in X$, $\mu(x+B)=0$. More generally, any set
$B\subset X$ is said Haar-nul if it is contained in a Haar-nul Borel set.
Let $U$ be an open subset of $X$. A set $P\subset U$ is said \emph{prevalent} in
$U$ if $U\setminus P$ is a Haar-nul set of $X$. To our knowledge, the first study of prevalence as a notion of ``large sets'' of Banach spaces goes back to \cite{HSY}. See \cite{Ott-Yorke} for a review on prevalence. As we can note, if $X$ is finite-dimensional, then a prevalent set is exactly a set of full Lebesgue measure. This means that a generic set may be negligible for the point of view of prevalence.

Of course, having different conclusions about the importance of a subset depending on the point of view is troublesome. However, we claim that the generic sets of the present article are also prevalent: they are ``large'' for both notions. Indeed, our proofs use Henry's Theorem (Theorem \ref{th_Henry}) and it is proved in \cite{RJ-prevalence} that this type of transversality theorem can be adapted to prevalence: the generic set obtained in its conclusion is also a prevalent set. 

\subsection{Allowing manifolds as state spaces}\label{sect_discuss_manifolds}

In the present article, we assume that the cell state spaces $X_c$ are of the type $\RR^{d_c}$. It is possible to extend the results to the case where $X_c$ are $\Cc^k-$manifolds of dimension $d_c$. Indeed, we can see that the central arguments of our proofs are purely local: we simply perturb the vector field along a small part of a trajectory. So, most of our arguments can be adapted to the case of manifolds simply by translating them to local charts, see \cite{RJ} for results proved in tori and for a short detailed example in the case of general manifolds. Note that our proofs require to view spaces as $\Cc^0([0,1],X)$ as Banach manifolds. To this end, we can define a local chart along a curve $x(\cdot)$ as follows. First cover the image of $x(\cdot)$ by a finite number of charts $\Oc_i$ of $X$, each associated to a time interval $[\sigma_i,\tau_i]$. It is now sufficient to explain what is a neighborhood of $t\in[\sigma_i,\tau_i]\mapsto x(t)$ in $\Cc^0([\sigma_i,\tau_i],\Oc_i)$. This is done by pulling back $\Cc^0([\sigma_i,\tau_i],\Oc_i)$ to $\Cc^0([\sigma_i,\tau_i],\RR^{d})$ through the chart. As one can see, working with manifolds only means heavier framework and notations. But this should not bring any obstruction to our main arguments, so we claim that our results extend to the case where the state spaces are general manifolds.

\subsection{About the self-dependence}\label{sect_discuss_self_depend}

In the present article, we choose to keep the framework of the previous works. 
In particular, we assume that the state $x_c$ of a cell $c$ is always a distinguished input of itself since $f_c$ is of the form $\hat f_c(x_c,x_{T(I(c))})$. This can be expressed as a particular case of more general networks with $f_c=\hat f_c(x_{T(I(c))})$. Making the self-dependence not automatic opens the possibility of $c$ not being an input of itself or being an input of itself but through an arrow associated with other ones.

In \cite{RJ,Maxime}, we obtain parts of our results without assuming the self-dependence. However, it was mandatory for several results, in particular when considering equilibrium points. It is thus possible that Proposition \ref{prop_synchrony_eq} fails without the self-dependence structure. For the sake of simplicity, and because it is the classical choice of the previous works, we choose to assume the self-dependence in all the present article. 
As already discussed in \cite{RJ} there are other arguments in favor of this choice:
\begin{enumerate}[1)]
\item If we allow the state spaces $X_c$ to be manifolds, then the vector fields depend of the position (since in general, the tangent space depends on the position). As discussed above, considering manifolds is simply heavier but all the results of the present article should remain true. This generalization cannot be done without self-dependence.
\item In biological models, it is difficult to imagine a cell or an individual, for which the evolution of its state is independent of itself. Even in models where this self-dependence is very light, small perturbations of the vector field with respect to this self-dependence are reasonable. Remember that cells of the same type have the same self-dependence, so we do not destroy the structure and symmetries of the networks by doing so.
\item If there is a scientific field where the self-dependence in ODEs is not automatic, it is clearly physics. Therefore, the self-dependence may be irrelevant in ODEs defined by physical laws. However, in this case, the entire vector field is itself highly constrained, so proving that a result is generic with respect to the vector field is certainly irrelevant. To study a physical model, we must precisely define its structure and the parameters that may vary (masses, distances\ldots but probably not the power of the law) and prove an ad hoc result in this specific class.   
\end{enumerate}


\end{document}